\documentclass[11pt]{article}
\usepackage{amsmath}
\usepackage{amsfonts}
\usepackage{amssymb}
\usepackage{latexsym, amscd, amsbsy, xypic, color,verbatim}
\usepackage{mathrsfs,bbm,dsfont}
\usepackage{graphicx, graphics, pstricks}
\usepackage{caption}
\usepackage[mathscr]{eucal}
\usepackage[all]{xy}
\usepackage{fancyhdr}
\usepackage{fullpage}
\usepackage{setspace}
\usepackage{subfigure}
\def\thefootnote{\fnsymbol{footnote}}
\newtheorem{thm}{Theorem}[section]
{}
\newtheorem{prop}[thm]{Proposition}
\newtheorem{lemma}[thm]{Lemma}

\newtheorem{defin}[thm]{Definition}
\newtheorem{remark}{remark}[section]

\newcommand{\proof}{{\it Proof. \quad}}
\newcommand{\qed}{\hfill\Box\medskip}

\usepackage{CJK}

  \def\De{\Delta}

 \def\de{\delta}

\def\bbn{{\mathbb N}}

  \def\leq{\leqslant}  \def\geq{\geqslant}

\def\lmto{\longmapsto}

\def\Hom{\mbox{\rm Hom}}  
 \def\id{\mbox{\rm id}}
 \def\res{\mbox{\rm Res}}
  \def\ind{\mbox{\rm Ind}}
    \def\supp{\mbox{\rm supp}\,}
  
\def\Ext{\mbox{\rm Ext}\,}   
  
\def\dim{\mbox{\rm dim}\,}   \def\End{\mbox{\rm End}\,}
 
\def\mod{\mbox{\rm {\textbf{-mod}}}\,}

\def\rad{\mbox{\rm rad}\,}   \def\red{{\text{\rm red}}}

\def\ggp#1#2{\left[\kern-3.2pt\left[{#1\atop #2}\right]\kern-3.2pt\right]}

\def\A{\mathcal{A}}

\def\la{\lambda}
\def\La{\Lambda}

\def\de{\delta}

\def\al{\alpha}
\def\pr{\prime}

\def\o{\otimes}
\def\bo{\bigotimes}

\begin{document}

\title{\bf Induction and restriction functors for cellular categories}

\author{
Pei Wang }
\date{}

 \maketitle

\begin{abstract}
 Cellular categories are a generalization of cellular algebras, which include a number of important categories
such as 
(affine)Temperley-Lieb categories, Brauer diagram categories, partition categories, the categories of invariant tensors for certain quantised enveloping algebras and their highest weight representations, Hecke categories and so on. The common feather is that, for most of the examples, the endomorphism algebras of the categories  form a tower of algebras. In this paper, we give an axiomatic framework for the cellular categories related to the  quasi-hereditary tower and then study the representations in terms of induction and restriction. In particular, a criteria for the semi-simplicity of  cellular categories is given by using the cohomology groups of cell modules.
 Moreover, we investigate the algebraic structures on Grothendieck groups of cellular categories and provide a diagrammatic approach to compute
the multiplication in the Grothendieck groups of Temperley-Lieb categories.

\medskip

\end{abstract}

\renewcommand{\thefootnote}{\alph{footnote}}
\setcounter{footnote}{-1} 
\renewcommand{\thefootnote}{\alph{footnote}}
\setcounter{footnote}{-1} 

\section{Introduction}
Cellular categories were  defined by Westbury in \cite{Wes3} as a generalization of cellular algebras, which were first introduced by Graham and Lehrer in \cite{GL}.
In a cellular category, the $\hom$-space of any two objects is spanned by a distinguished basis, so-called cellular basis. Therefore, an endomorphism algebra in a cellular category is cellular. In particular, if we can regard an algebra as a category with one object, then a cellular algebra is indeed a cellular category.
It was shown that many important classes of algebras arising in representation theory, invariant theory, knot theory, subfactors and statistical mechanics are cellular (see e.g. \cite{Bow,Ge,GL,KC-1,Mur,Ng,XiB,XiP}), and  most of their categorical analogues  are also cellular, such as 
Temperley-Lieb categories \cite{Wes2}, Brauer diagram categories \cite{LZ}, partition categories \cite{Jo,Mar3}, the categories of invariant tensors for certain quantised enveloping algebras and their highest weight representations \cite{Wes3,Wes1}, the categories of Soegrel bimodules \cite{EW} and other more general Hecke categories (a strictly object-adapted cellular category due to \cite{EL}).
As is known,
if an algebra admits a cellular structure, one will have a practicable way to describe the  representations and homological properties of the algebra \cite{GL,CY}.
In this artical,  we shall investigate the properties of cellular categories.

Lots of important examples of cellular algebras actually occur in towers $A_{0}\subset A_{1}\subset A_{2} \subset \cdots$ with an intense interplay between each other in terms of induction and restriction. 
These include  Temperley-Lieb algebras \cite{MartinB} and their cyclotomic analogues \cite{RX}, Brauer algebras \cite{Bra,Bro,Mk,Mk2}, blob algebras \cite{Ms,MW}, partition algebras \cite{Bl,G0,Mar3,MarE} and so on. 
The tower method  was first developed by Jones \cite{Jo0} and Wenzl \cite{Wen} for semi-simple case. Further, for the general case, Cox, Martin,
Parker and Xi established a framework of towers of quasi-hereditary algebras by combining the ideas from the tower formalism in \cite{GHJ} with the notion of recollement in \cite{CPS}. Then, influenced by the work of K\"{o}nig and Xi \cite{KC0} as well as by the work of Cox et al\cite{CMPX}, the analogous for cellularity were given by Goodman and Graber\cite{GG}.

There are in fact a number of cellular categories with the endomorphism algebras forming a tower. Further, note that the $\hom$-space of two objects in a cellular category is a natural nexus for their endomorphism algebras.

Those motivate us to introduce a class of cellular categories by an axiomatic menner, so-called \textit{cellular tower category}, or \textit{CTC} (see Section 3.2), in which the endomorphism algebras are required to be quasi-hereditary,   and herein the  induction and restriction behave well in the tower. In this setting, the homological aspects of representation theory are computed efficiently by using induction and restriction.

Precisely, let $K$ be a field, and $\A$~be a cellular tower category (see Definition \ref{deftower}) with its  object set  the set of natural numbers  $\mathbb{N}$. By definition, the endomorphism algebra $A_{n}$ of an object $n$ is a cellular algebra with an index set
$\La_{n}$. Suppose that $\De_{n}(\la)$ denotes the cell module of $A_{n}$ corresponding to the index $\la \in \La_{n}$.

Our main result  can be stated as the following theorem.

\begin{thm}\label{semis1}
Let~$\A$~be a cellular tower category.  Suppose that for all~$n \in \bbn$~and pairs of indices
~$\la \in \La_{n}\backslash \La_{n-2} ~and ~ \mu \in \La_{n}\backslash \La_{n-4}$~we have$$\Ext_{\scriptsize{A_{n}}}^{1}(\De_{n}(\la),\De_{n}(\mu) )=0,$$
Then each of the endomorphism algebras~$A_{n}$~in $\A$~is semi-simple.
\end{thm}

In \cite{CY}, Cao provided a criteria of semi-simplicity for a  cellular algebra by checking  the first cohomology groups of cell modules
for all indices. Rather,
Theorem 1.1 tells us that, the verification needs only  some of the indices in a cellular tower category.

The Grothendieck groups of the categories of finitely generated modules and finitely generated projective modules over a tower of algebras can be endowed with algebra and coalgebras structures.
Many cases of interest, such as symmetric group algebras, Hecke algebras and other deformations,  give rise to a dual pair of Hopf algebras, which are realization of some classical algebras in the theory of symmetric function \cite{Zsky,KT,BHT,Se}. The common feature is that the examples admit Mackey's formula, which implies that the comultiplication is an algebra homomorphism, that is, making both Grothendieck groups into bialgebra.

In \cite{BL}, Bergeron and Li gave an analogue of Mackey's formula by an axiom and introduced a general notation of a tower of algebras, which ensures  that the Grothendieck groups of a tower of algebras can be a pair of graded dual Hopf algebras. For a cellular category with a quasi-hereditary tower, we shall see that the analogue of Mackey's formula never holds. The Grothendieck groups, however, have algebra and coalgebra structures under  certain conditions. Further, we shall study the algebraic structures on the Grothendieck groups of Temperley-Lieb categories.

This paper is organized as follows. In Section 2, we shall recall some notations and basic facts.
In Section 3, we first study  Morita contexts of endomorphism algebras in  a cellular category in 3.1 and then  we introduce the cellular tower categories and prove Theorem \ref{semis1} in 3.2.
Section 4 studies the algebraic structures on Grothendieck groups of cellular categories. In particular, we give a method to compute the multiplication in the Grothendieck group of a Temperley-Lieb category in 4.2.

\section{Preliminaries}
In this section, we shall recall some basic definitions and facts needed in our later proofs.

Throughout the paper,  all algebras are finite-dimensional algebras over a fixed field $K$. All modules are finitely generated unitary left modules. For an algebra A, the category of A-modules is denoted by $A\mod$.
Let $\A$  always be a small $K$-linear category with finite dimensional hom-spaces, that is, the class of objects is a set and every hom-set is a  finite dimensional $K$-vector space and the composition map of morphisms is bilinear. 
For an object $n$
in $\A$, the endomorphism algebra $\End (n)$ of $n$ is
simply denoted by $A_{n}$.

\subsection{Cellular categories }

We now recall the definition of cellular categories, which are generalised by Westburry \cite{Wes3} from 
the cellular algebras.

\begin{defin}{\textsc{\cite{Wes3}}}\label{basis}
 Let~$\mathcal{A}$~be a~$K$-linear category with an anti-involution
~$*$~( that is,~$*$~is a dual functor on~$\mathcal{A}$~with~$(-)^{**}=\id_{\A}$). 
Then cell datum for~$\mathcal{A}$~consists of a 
 partially ordered set~$\Lambda$,   a finite set~$M(n,\lambda)$ for each~$\lambda \in\Lambda$~and each object~$n$~of~$\mathcal{A}$, and for~$\lambda \in\Lambda$ and~$m, n$~any two objects~$\mathcal{A}$~we have an inclusion
$$C:  M(m, \lambda)\times  M(n, \lambda)\rightarrow\Hom_{\mathcal{A}}(m,n)$$
$$C: (S, T)\lmto C_{S, T}^{\lambda}.$$

The conditions for this datum are required to:

$(C1)$ For all objects~$m, n$~in~$\mathcal{A}$, the image of the map
$$C: \coprod_{\lambda \in\Lambda}M(m, \lambda)\times  M(n, \lambda)\rightarrow\Hom_{\mathcal{A}}(m,n)$$
is a basis for~$\Hom_{\mathcal{A}}(m,n)$~as a~$K$-space.

$(C2)$ For all objects~$m, n$, all~$\lambda \in\Lambda$~and~$S\in M(m,\lambda)$, $T\in M(n,\lambda)$~we have
$$(C_{S,T}^{\lambda})^{*}=C_{T,S}^{\lambda}.$$

$(C3)$ For all objects~$p, m, n$, all~$\lambda \in\Lambda$~and all~$a\in \Hom_{\mathcal{A}}(p,m)$,
$S\in M(m,\lambda)$,~$T\in M(n,\lambda)$~we have
$$aC_{S,T}^{\lambda} \equiv \sum_{S^{\prime}\in M(p, \lambda)} r_{a}(S, S^{\prime})C_{S^{\prime},T}^{\lambda} ~~\mbox{\rm mod}~\mathcal{A}(<\lambda),$$
where~$r_{a}(S^{\prime}, S)\in K$~is independent of~$T$~and~$\mathcal{A}(<\lambda)$~is the~$K$-span of
$$\{ C_{S,T}^{\mu} ~|~\mu<\lambda;~ S\in M(p,\mu), T\in M(n,\mu) \} .$$
\end{defin}


\textbf{Remark 2.1}
	If we regard any $K$-algebra as a $K$-linear category with one object, then this is a generalisation of the definition of a cellular algebra. On the other hand, suppose~$\A$~is a cellular category. For any object~$n$, let ~$\La_{n}=\{ ~\la\in \La\mid M(n,\la)\neq \emptyset~ \}$~and~$M_{n}:=\bigcup\limits_{\la \in \Lambda_{n}}M(n,\la)$, then the endomorphism algebra~$A_{n}$ of $n$~is a cellular algebra with cell datum~$(\Lambda_{n}, M_{n}, C, *)$, where ~$C$~and~$*$~are  restrictions on~$A_{n}$.
The basis $\{ C_{S, T}^{\lambda} ~|~ S, T \in M_{n},  \lambda \in \La_{n}   \}$~is called  a cellular basis for $A_{n}$.

In~\cite{GL}, Graham~and~Lehrer~introduced the definition of cell modules of a cellular algebra by using the cellular basis. For cellular categories,  we can define   cell modules of  the endomorphism algebra of an object  in the same way.

\begin{defin}
Let~$\A$~be a cellular category.  For each object~$n$, let~$\La_{n}$~be the index set of the cellular algebra~$A_{n}$~(see Remark 2.1 above ). For each~$\la\in \La_{n}$~define the left ~$A_{n}$-module~$\De_{n}(\la)$~as follows: $\De_{n}(\la)$~is a~$K$-space with basis ~$C^{(n,\la)}=\{ C_{X}^{(n,\la)}\mid X\in M(n,\la)\}$ and $A_{n}$-action defined by 
$$ aC_{X}^{(n,\la)} = \sum_{Y\in M(n, \lambda)} r_{a}(X, Y)C_{Y}^{(n,\la)},  ~~~~~~(a\in A_{n}, ~X\in M(n, \lambda))$$
where~$r_{a}(X, Y)$~is the element of $K$~defined in~$(C3)$.
$\De_{n}(\la)$~is called the cell module of ~$A_{n}$~corresponding to~$\la$.
\end{defin}

Applying `$*$' on ~$(C3)$ in Definition~\ref{basis}, we obtain
$$C_{T,S}^{\lambda}a^{*} \equiv \sum_{S^{\prime}\in M(p, \lambda)} r_{a}(S, S^{\prime})C_{T,S^{\prime}}^{\lambda} ~~\mbox{\rm mod}~\mathcal{A}(<\lambda). \leqno{(C3^{\pr})}$$
where~$r_{a}(S^{\prime}, S)\in K$~is independent of~$T$.

The following lemma is a direct generalization of ~\cite[Lemma 1.7]{GL} from cellular algebras to cellular categories.
\begin{lemma}\label{shang}
Let~$\A$~be a cellular category, $m,n~and~q$~be objects in $\A$, and let~$\la\in \La$. Then for any elements $U\in M(m,\la)$, $T,X \in M(n,\la)$ and $Y\in M(p,\la)$, we have
$$C_{U,T}^{\lambda} C_{X,Y}^{\lambda}\equiv\phi_{(n,\la)}(T,X) C_{U,Y}^{\lambda} ~~\mbox{\rm mod}~\mathcal{A}(<\lambda), $$
where $\phi_{(n,\la)}$ is a map from $M(n,\la)\times M(n,\la)$ to $K$.
\end{lemma}
\proof By~$(C3)$, we have
$$C_{U,T}^{\lambda} C_{X,Y}^{\lambda}  \equiv\sum_{Z\in M(m, \lambda)} r_{C_{U,T}^{\lambda}}(X, Z)C_{Z,Y}^{\lambda} ~~\mbox{\rm mod}~\mathcal{A}(<\lambda),$$
and by~$(C3^{\pr})$ , it follows
$$C_{U,T}^{\lambda} C_{X,Y}^{\lambda}  \equiv\sum_{V\in M(p, \lambda)} r_{C_{Y,X}^{\lambda}}(V, T)C_{U,V}^{\lambda} ~~\mbox{\rm mod}~\mathcal{A}(<\lambda).$$

Comparing  the previous equations follows
$$r_{C_{U,T}^{\lambda}}(X, U)C_{U,Y}^{\lambda}=r_{C_{Y,X}^{\lambda}}(Y, T)C_{U,Y}^{\lambda}.$$
Because~$r_{C_{U,T}^{\lambda}}(X, U)$~is  independ of~$Y$~by $(C3)$~and~$r_{C_{Y,X}^{\lambda}}(Y, T)$~is independ of~$U$ by $(C3')$,  we may write $$\phi_{(n,\la)}(T,X):=r_{C_{U,T}^{\lambda}}(X, U)=r_{C_{Y,X}^{\lambda}}(Y, T).$$ Hence $$C_{U,T}^{\lambda} C_{X,Y}^{\lambda}\equiv\phi_{(n,\la)}(T,X) C_{U,Y}^{\lambda} ~~\mbox{\rm mod}~\mathcal{A}(<\lambda). $$
$\qed$

Thus we  can define  a bilinear form~$\phi_{(n,\la)}:\De_{n}(\la)\times \De_{n}(\la)\rightarrow K$ by $$\phi_{(n,\la)}(C_{U}^{(n,\la)},C_{V}^{(n,\la)})=\phi_{(n,\la)}(U,V),$$
where~$C_{U}^{(n,\la)},C_{V}^{(n,\la)}\in C^{(n,\la)}$ with $U,V \in M(n,\la)$, extended $\phi_{(n,\la)}$ bilinearly.

The following lemma collects some known facts for the bilinear form $\phi_{(n,\la)}$.
\begin{lemma}\label{xing}~\cite[Prop 2.4]{GL}
Keep the notation above. Then$:$

$(1)$~$\phi_{(n,\la)}$~is symmetric, that is,~$x,y \in \De_{n}(\la)$, $\phi_{(n,\la)}(x,y)=\phi_{(n,\la)}(y,x)$$;$

$(2)$~For~$x,y \in \De_{n}(\la)$~and~$a\in A_{ n}$, we have~
$\phi_{(n,\la)}(a^{*}x,y)=\phi_{(n,\la)}(x,ay)$$;$

$(3)$~For~$C_{S,T}^{\la}\in A_{n}$~and~$C_{U}^{(n,\la)}\in C^{(n,\la)}$, we have
$C_{S,T}^{\la}\cdot C_{U}^{(n,\la)}=\phi_{(n,\la)}(T,U)C_{S}^{(n,\la)}$.

\end{lemma}

Define
$$\rad_{n}(\la):=\{x\in \De_{n}(\la)~|~\phi_{(n,\la)}(x,y)=0 ~for~all~y\in \De_{n}(\la)\}.$$
If $\phi_{(n,\la)}\neq 0$, then $\rad_{n}(\la)$ is the radical of the $A_{n}$-module $\De_{n}(\la)$.

Let $\La_{n}^{0}=\{\la \in \La_{n}~|~ \phi_{(n,\la)}\neq 0\}$. The following result shows that this set parameterizes the simple modules, which was proved by Graham and Lehrer for cellular algebras.

\begin{thm}\cite[Prop 3.4]{GL}
Let $A$ be a cellular $K$-algebra with the cell datum $(\La, M, C, *)$. Suppose $\De(\la)$ and $ \phi_{\la}$ are the cell module and  the bilinear form, respectively, corresponding to $\la \in \La$. Let $\La^{0}=\{\la \in \La~|~ \phi_{\la}\neq 0\}$.
Then the set $\{ L(\la):= \De(\la)/ \red(\la)~|~ \phi_{\la}\neq 0\}$ is a complete set of non-isomorphic absolutely simple $A$-module.

\end{thm}

The following theorem says that the issue of semi-simplicity reduces  to the computation of the discriminants of bilinear forms associated to cell modules. 

\begin{thm}\label{sem}\cite[Prop 3.8]{GL}
Let $A$ be a cellular $K$-algebra as above. Then the following are equivalent:

(1) The algebra $A$ is semi-simple;

(2) All cell modules are simple and pairwise non-isomorphic;

(3) The bilinear form $\phi_{\la}$ is non-degenerate (that is, $\rad(\la)=0$) for each $\la \in \La$
\end{thm}

In \cite{KC7,KC8}, Konig and Xi investigated the relationships between cellular algebras
and quasi-hereditary algebras in terms of comparing the so- called cell chains with hereditary chains (for   quasi-hereditary algebras we refer to \cite{CPS}), they gave some  criterions for a cellular algebra to be quasi-hereditary. 

The following is equivalent to \cite[Theorem 3.1]{KC7}.

\begin{thm} \label{shang}
Let~$A$~be a cellular algebra with the cell datum $(\La, M, C, *)$, and $\La^{0}  $ be as above. Then $A$ is quasi-heredity  if and only if    $\La=\La^{0}$.
\end{thm}

This lemma says that a cellular algebra
 is quasi-heredity  if and only if   the poset $\La$ coincides with $\La^{0}$, that is,  $\phi_{\la}\neq 0$ for all $\la\in\La$. Thus in the case the set $\La$ parameterizes the simple modules.

For our purpose in this paper, we need that the endomorphism algebras of a cellular category are quasi-hereditary, so we define:

\begin{defin}\label{hcc}
Let~$\A$~be a cellular category.   Suppose that

(1)  the object set of $\A$ is the set of natural numbers  $\mathbb{N}$,

(2) for all $m ,n \in \mathbb{N}$ satisfying  $m \leqslant n$, if $\Hom_{\A}(m,n)\neq 0$~(equivalently, $\Hom_{\A}(n,m)\neq 0$ by anti-involution $*$), then $\La_{m}$ is a saturated ordered subset of $\La_{n}$, that is, $\La_{m}\subset \La_{n}$ preserving ordering, and if $\la < \mu$ with~ $\mu\in\La_{m}$, $\la\in\La_{n}$, then $\la\in\La_{m}$.

(3) each endomorphism algebra $A_{n}$ is quasi-hereditary.

Then~$\A$~is called a hereditary cellular category.
\end{defin}

\textbf{Remark 2.2}
(A)~~In Definition \ref{hcc}, condition (2) says that index sets $ \La_{n}$ preserve the ordering of natural numbers.

(B)~~Let~$\A$~be a hereditary cellular category with cell datum~$(\Lambda, M, C, *)$.
For all~$m, n\in \mathbb{N}$, denote by~$\La_{(m,n)}:=\{ ~\la\in \La\mid M(m,\la)\neq \emptyset~and~M(n,\la)\neq \emptyset\}$, and thus $\La_{(n,n)}$ 
is just $\La_{n}$.

By condition (2), if $\Hom_{\A}(m,n)\neq 0~with~m \leqslant n$ and~$M(m,\la)\neq \emptyset$, then~$\la \in \La_{m}\subseteq \La_{n}$. Therefore this imples~$M(n,\la)\neq \emptyset$ and~$\La_{(m,n)}=\La_{m}$.

\section{Induction and restriction functors for cellular categories}
Let $\A$ be a $K$-linear category, and let $m,n$ be objects in $\A$. As is known, the hom-space $\Hom(m,n)$
is   a natural left $A_{ m}$-right $A_{n}$-bimodule, and hence we can consider the functors $\Hom_{\A}(m,n)\o_{\scriptsize A_{ n}}-$ and $-\o_{\scriptsize A_{  m}}\Hom_{\A}(m,n)$.  In this section, we first study those functors in Section 3.1.

In Section 3.2, We first give the definition of cellular tower category by an axiomatic menner and then give a proof of our main result Theorem 1.1.

\subsection{Morita contexts}

Let~$\A$~be a hereditary cellular category with cell datum~$(\Lambda, M, C, *)$.
For each~$m, n\in \mathbb{N}$, $\Hom_{\A}(m,n)\o_{\scriptsize A_{n}}-$~is a functor from category~$A_{n}\mod$~to category~$A_{m}\mod$. For each $\la\in\La_{(m,n)}=\{ ~\la\in \La\mid M(m,\la)\neq \emptyset~and~M(n,\la)\neq \emptyset\}$, let~$\De_{m}(\la)$ and $ \De_{n}(\la)$~be the cell modules of~$A_{m}$~and
$A_{n}$ corresponding to~$\la$,~respectively. 

It is easy to see that the following map is an $A_{m}$-module homomorphism$:$
$$\alpha:  \Hom_{\A}(m,n)\o_{\scriptsize A_{n}} \De_{n}(\la)\rightarrow \De_{m}(\la)$$
$$~~~~~~~~~~~~~~~~~ a \o C_{X}^{(n,\la)}~~~~\mapsto ~~~~\sum\limits_{Z\in M(m,\la)}r_{a}(Z,X) C_{Z}^{(m,\la)}$$
where~$a\in \Hom_{\A}(m,n), C_{X}^{(n,\la)}\in C^{(n,\la)}$,~$r_{a}(Z,X)$~ is given by definition \ref{basis}(C3).

The following lemma show that it is further an $A_{m}$-module isomorphism.

\begin{lemma}\label{hom}\label{hanzi1}
Let~$\A$~be a hereditary cellular category. For each~$m, n\in \mathbb{N}$~and~$\la\in \La_{(m,n)}$,
we have$$\Hom_{\A}(m,n)\o_{\scriptsize A_{n}} \De_{n}(\la)\simeq \De_{m}(\la),$$
as an~$A_{m}$-isomorphism.
\end{lemma}
\proof
Because~$A_{n}$~is quasi-hereditary,  it follows~$\phi_{(n,\la)}\neq 0$ by Theorem \ref{shang}, hence there exist~$U_{0}, T_{0}\in M(n,\la)$~such that~$\phi_{(n,\la)}(U_{0}, T_{0})\neq 0$. We then fix $U_{0}, T_{0}$.

For any~$C_{Y}^{(m,\la)}\in C^{(m,\la)}$~with~$Y\in M(m,\la)$,
we have~$\frac{1}{\phi_{(n,\la)}(U_{0},T_{0})} \al(C_{Y,U_{0}}^{\la}\o C_{T_{0}}^{(n,\la)})=C_{Y}^{(m,\la)}$. It follows that~$\al$~is surjective and 
$$\dim_{K} (\Hom_{\A}(m,n)\o_{\scriptsize A_{n}} \De_{n}(\la)) \geq \dim_{K} (\De_{m}(\la))= ~~^\# M(m,\la).$$

To prove previous inequality is actually  a equality, fixed $U_{0}, T_{0}$ as above, it is sufficient to show that~$\{ C_{S,U_{0}}^{\la}\o C_{T_{0}}^{(n,\la)}\mid S\in M(m,\la)\}$ is a spanning set~of the~$K$-space~$\Hom_{\A}(m,n)\o_{\scriptsize A_{n}} \De_{n}(\la)$. Thus, for any~$a\in \Hom_{\A}(m,n)$~and~$C_{X}^{(n,\la)}\in C^{(n,\la)}$, we have
\begin{eqnarray*}
a\o C_{X}^{(n,\la)} & =&\frac{1}{\phi_{(n,\la)}(U_{0},T_{0})}
~a\o C_{X,U_{0}}^{\la}\cdot C_{T_{0}}^{(n,\la)} \\
&=& \frac{1}{\phi_{(n,\la)}(U_{0},T_{0})}
~a C_{X,U_{0}}^{\la}\o C_{T_{0}}^{(n,\la)}\\
&=& \frac{1}{\phi_{(n,\la)}(U_{0},T_{0})}(
~\sum\limits_{S\in M(m,\la)}r_{a}(X,S) C_{S,U_{0}}^{\la}\o C_{T_{0}}^{(n,\la)}+ b \o C_{T_{0}}^{(n,\la)}),
\end{eqnarray*}
where~$b\in \mathcal{A}(<\lambda)$. We next show~$ b \o C_{T_{0}}^{(n,\la)}=0$.

Without loss of generality, consider a basis element~$b=C_{W,V}^{\mu}$~with~$\mu < \la,~W\in M(m,\mu)~and~V\in M(n,\mu)$. According to the quasi-heredity, we have~$\phi_{(n,\mu)}\neq 0$. Hence, there exist~$U^{\pr}, T^{\pr}\in M(n,\mu)$ such that~$\phi_{(n,\mu)}(U^{\pr}, T^{\pr})\neq 0$. It follows
\begin{eqnarray*}
C_{W,V}^{\mu}\o C_{T_{0}}^{(n,\la)} & =&\frac{1}{\phi_{(n,\mu)}(U^{\pr},T^{\pr})}
~C_{W,U^{\pr}}^{\mu} C_{T^{\pr},V}^{\mu}\o C_{T_{0}}^{(n,\la)} \\
&=& \frac{1}{\phi_{(n,\mu)}(U^{\pr},T^{\pr})}
~C_{W,U^{\pr}}^{\mu}\o C_{T^{\pr},V}^{\mu}\cdot C_{T_{0}}^{(n,\la)}=0.
\end{eqnarray*}

Hence $\alpha$ is injective. This finishes the proof.
$\qed$

This lemma says that~the functor $\Hom_{\A}(m,n)\o_{\scriptsize A_{n}}-$~sends a cellular module to a cellular module with the same index. In fact, we shall show that this functor is an idempotent embedding
functors \cite{ASS}. We first recall the definition of~Morita context.

\begin{defin}~{\rm\cite{Lam}}
Let~$A, B$~be two~$K$-algebras, and let $_B P_A$,~$ _A Q_B$~be bimodules, and $\theta, \phi$ be a pair of bimodule homomorphisms
$$\theta : P\o_{A}Q\rightarrow  _{B}B_B,~~~~~~~\phi : Q\o_{B}P\rightarrow  _{A}A_A$$
such that for all~$x,y\in P$~and~$f,g\in Q$, 
$$\theta(x\o f)y=x\phi(f\o y),~~~~~~f\theta(x\o g)=\phi(f\o x)g.$$

Then the tuple~$(A, _B P_A, _A Q_B, B, \theta, \phi)$~is called a~Morita context.
\end{defin}

If~$\theta$~is surjective in  the previous definition, we have the following lemma.

\begin{lemma}~{\rm\cite[Prop 18.17]{Lam}}~\label{morita} Keep the notation above. 

Suppose~$\theta$~is surjective. Then

$(1)$~$\theta$~is an isomorphism$;$

$(2)$~$P_A$~and~$_AQ$~are finitely generated projective.

$(3)$~There are $K$-algebra isomorphism~$B\simeq\End(P_A)\simeq\End(_AQ)$.
\end{lemma}

Let~$\A$~be a hereditary cellular category. For objects~$m,n\in\mathbb{N}$,  the tuple~$$(A_{m},~ _{\scriptsize{A_{m}}}\Hom_{\A}(m,n)_{\scriptsize{A_{n}}}~, ~_{\scriptsize{A_{n}}}\Hom_{\A}(n,m)_{\scriptsize{A_{m}}}~, A_{m} ,~ \eta,~\rho)$$ is a~Morita context, where~$\Hom(m,n)$~is a $A_{m}$-$A_{n}$-bimodule, ~$\Hom(n,m)$ is a $A_{n}$-$A_{m}$-module. The map~$$\eta :~  \Hom_{\A}(m,n)\o_{\tiny{A_{n}}}\Hom_{\A}(n,m)\rightarrow A_{m}$$
defined by the composition of morphisms, is an~$A_{m}$-bimodule homomorphism;

Similarly, the map~$$\rho:~  \Hom_{\A}(n,m)\o_{\tiny{A_{m}}}\Hom_{\A}(m,n)\rightarrow A_{n}$$
defined by the composition of morphisms, is an~$A_{n}$-bimodule homomorphism.

\begin{lemma}~\label{sur} Keep the notation above. 
Let~$m,n\in\mathbb{N}$ 
satisfying  $m \leqslant n$. Suppose $\Hom_{\A}(m,n)\neq 0$.
Then~$\rho$~is surjective.

\end{lemma}
\proof Let $C_{S,T}^{\la}$ be a cellular basis element of~$A_{n}$ with $\la \in \La_{n}$ and $S, T\in M(n,\la)$.
By Definition \ref{hcc} and Remark 2.2(B), we have  $\La_{n}=\La_{(n,m)}\subseteq \La_{m}$. Hence $\la \in \La_{m}$. Since~$A_{m}$~is quasi-hereditary, it follows~$\phi_{(m,\la)}\neq 0$ by Theorem \ref{shang}, therefore there exist~$U_{0}, V_{0}\in M(m,\la)$ such that~$\phi_{(m,\la)}(U_{0}, V_{0})\neq 0$.
This implies~$$\rho (\frac{1}{\phi_{(n,\la)}(U_{0},V_{0})} C_{S,U_{0}}^{\la}\o
C_{V_{0},T}^{\la})= C_{S,T}^{\la},$$
Because~$C_{S,T}^{\la}$~is arbitrary, $\rho$~is surjective.

$\qed$

As an immediate consequence of lemma \ref{morita} and \ref{sur} we get the following.

\begin{prop}\label{hanzi2}
Let~$\A$~be a hereditary cellular category. Suppose~$m,n\in\mathbb{N}$~satisfying~$m\leq n$~and~$\Hom_{\A}(m,n)\neq 0$. Then

$\mbox{\rm (}1\mbox{\rm )}$~~$\rho$~is an~$A_{n}$-bimodule isomorphism$;$

$\mbox{\rm (}2\mbox{\rm )}$~~There exist  idempotents $e$~and~$f$~of~$A_{m}$~such that
~$\Hom_{\A}(m,n)\simeq(A_{m}) e$~as left~$
A_{m}$-module isomorphism, and~$\Hom_{\A}(n,m)\simeq f(A_{m}) $~as right~$
A_{m}$-module isomorphism, that is,~$\Hom_{\A}(m,n)\o_{\tiny{A_{  n}}} -$~and~$-\o_{\tiny{A_{  n}}} \Hom_{\A}(n,m)$~are idempotent embedding functors;

$(3)$~~There are $K$-algebra isomorphism~$A_{n}\simeq e(A_{m})e\simeq f(A_{m})f$.
\end{prop}$\qed$

\subsection{Cellular tower categories}

Inspired by the theory of quasi-heredity towers of recollement studied  by Cox, Martin, Parker and Xi,
we first introduce the definition of  cellular tower category. Under the framework, the homological aspects of representation theory are computed efficiently by using induction and restriction. We then give a criteria for a cellular tower category to be semi-simple.

\begin{defin}\label{deftower}
	A hereditary cellular category is called a  cellular tower category (or CTC) if it satisfies~$\mbox{\rm (}A1\mbox{\rm )}$--~$\mbox{\rm (}A4\mbox{\rm )}$ as follow:
\end{defin}

$\mbox{\rm (}A1\mbox{\rm )}$ {\bf For each~$n \geq 0$~the endomorphism algebra~$A_{n}$~ can be identified with  a subalgebra of~$A_{ n+1}$~which preserves the identities.}

For an~$A_{ n+1}$-module~$M$, it has a natural ~$A_{n}$-module structure. Furthermore, we have the restriction functor:
$$\res_{n+1} : ~~ A_{ n+1}\mod\rightarrow A_{n}\mod$$
$$~~~~~~~~~~~~~~~~~_{\scriptsize A_{  n+1}}M~~~\mapsto~~~ _{\scriptsize A_{n}}M=\Hom_{\scriptsize A_{ n+1}}({A_{  n+1}} _{\scriptsize A_{n}},M).$$

We also have the induction functor:
$$\ind_{n} : ~~ A_{n}\mod\rightarrow A_{ n+1}\mod$$
$$~~~~~~~~~~~~~~~~~_{\scriptsize A_{n}}M~~~\mapsto~~~ A_{ n+1}\o_{\scriptsize A_{n}}M.$$

Our next three axioms ensure that induction and restriction behave well.

$\mbox{\rm (}A2\mbox{\rm )}$ {\bf For all~$n \geq 2$~we have the following $ A_{n-1}$-$ A_{n-2}$-bimodule isomorphism:}
$$\Hom_{\A}(n,n-2)\simeq_{\scriptsize A_{n-1}} {A_{n-1}}_{\scriptsize A_{n-2}}.$$

By our assumption that~$\A$~is quasi-hereditary, due to lemma~\ref{hanzi1}, for any~$\la \in \La_{n}$ we have
$\Hom_{\A}(n+2,n) \o \De_{n}(\la)\simeq \De_{n+2}(\la)$.  Furthermore, by using $\mbox{\rm (}A1{\rm )}$~and~$\mbox{\rm (}A2\mbox{\rm )}$, we have
\begin{equation}\label{1a}
\ind_{n}(\De_{n}(\la))= \res_{(n-1)}(\Hom_{\A}(n+2, n)\o_{\scriptsize A_{n}} \De_{n}(\la)).
\end{equation}

If an~$A_{n}$-module~$M$~in~$A_{n}\mod$ has a ~$\De_{n}$-filtration, that is, a filtration with successive quotients isomorphic to some  cell modules~$\De_{n}(\la)$'s, then we define the \textit{support} of~$M$, denoted by $\supp_{n}(M)$, to be the set of labels $\la$ for which $\De_{n}(\la)$ occurs in this filtration.

$\mbox{\rm (}A3\mbox{\rm )}$ {\bf For all~$m, n\in \bbn$~satisfying~$m\leq n$~and that~$n-m$~is even, and for all~$\la \in \La_{m}\backslash \La_{m-2}$, we have that~$\res_{n}(\De_{n}(\la))$~has a~$\De_{(n-1)}$-filtration and $$\supp_{n-1}(\res_{n}(\De_{n}(\la)))\subseteq (\La_{m-1}\backslash \La_{m-3} )\cup (\La_{m+1}\backslash \La_{m-1} )\subseteq  \La_{n-1}.$$}

For~$\mbox{\rm (}A3\mbox{\rm )}$, using~$\Hom_{\A}(n+2,n) \o \De_{n}(\la)\simeq \De_{n+2}(\la)$~and equation~(\ref{1a}), we deduce that for all~$m, n\in \bbn$~satisfying~$m\leq n$~and that~$n-m$~is even, and for all~$\la \in \La_{m}\backslash \La_{m-2}$, we have that~$\ind_{n}(\De_{n}(\la))$~has a~$\De_{(n+1)}$-filtration and$$\supp_{n+1}(\ind_{n}(\De_{n}(\la)))\subseteq (\La_{m-1}\backslash \La_{m-3} )\cup (\La_{m+1}\backslash \La_{m-1} )\subseteq  \La_{n+1}.$$

~$\mbox{\rm (}A4\mbox{\rm )}$ { \bf
Let~$n\in \bbn$. For each~$\la \in \La_{n}\backslash \La_{n-2}$~there exists~$\mu \in \La_{n-1}\backslash \La_{n-3}$~such that}
$$\la \in \supp_{n}(\ind_{n-1}(\De_{n-1}(\mu))).$$

By using~$\mbox{\rm (}A3\mbox{\rm )}$, $\mbox{\rm (}A4\mbox{\rm )}$~is equivalent to

~$\mbox{\rm (}A4^{\pr}\mbox{\rm )}$
Let~$n\in \bbn$. For each~$\la \in \La_{n}\backslash \La_{n-2}$, there exists~$\mu \in \La_{n-1}\backslash \La_{n-3}$~such that
$$\la \in \supp_{n}(\res_{n+1}(\De_{n+1}(\mu))).$$

For a quasi-hereditary algebra we have that~$\Ext^{1}(\De(\la),\De(\mu))\neq 0$~implies that~$\mu < \la$. Therefore~$\mbox{\rm (}A4\mbox{\rm )}$~is also equivalent to:

Let~$n\in \bbn$. For each~$\la \in \La_{n}$~there exists~$\mu \in \La_{n-1}$~such that there is a surjection
$$\ind_{n-1}(\De_{n-1}(\mu))\rightarrow \De_{n}(\la)\rightarrow 0.$$

\

Due to Cao~\cite{CY}, the following theorem provides some homological characterizations
of the semi-simplicity of cellular algebra.

\begin{thm}~\mbox{\rm{\cite[Thm.~1.2]{CY}}}\label{cao}
Let~$A$~be a cellular~$K$-algebra with respect to an  involution~$i$ and a poset~$(\La, \leq)$.
Let~$\La_{0}$~be the subset of~$\La$, which parametrizes the isomorphism classes of simple $A$-modules. 
Then the following statements are equivalent: 

~$\mbox{\rm (}a\mbox{\rm )}$~~~The algebra~$A$~is semi-simple.

~$\mbox{\rm (}b\mbox{\rm )}$~~~$\Ext_{A}^{1}(\De(\la), S(\mu))=0$~for any ~$\la, \mu \in\La_{0}$~satisfying ~$\mu \leq \la$, where~$S(\mu)$~ is the simple module with respect to ~$\mu$.

~$\mbox{\rm (}c\mbox{\rm )}$~~~$\Ext_{A}^{1}(\De(\la), \De(\mu))=0$~for any~$\la, \mu \in\La_{0}$~satisfying~$\mu \leq \la$.

~$\mbox{\rm (}c^{\pr}\mbox{\rm )}$~~\hspace{1pt}$\Ext_{A}^{1}(\De(\la), \De(\mu))=0$~for any~$\la, \mu \in\La$~satisfying~$\mu \leq \la$.

~$\mbox{\rm (}c^{\pr\pr}\mbox{\rm )}$~~$\Ext_{A}^{1}(\De(\la), \De(\mu))=0$~for all ~$\la, \mu \in\La$.

\end{thm}

The main result of this note can be stated as the following Theorem \ref{semis}.
\begin{thm}\label{semis}
Let~$\A$~be a cellular tower category.  Suppose that for all~$n \in \bbn$~and pairs of indices
~$\la \in \La_{n}\backslash \La_{n-2} ~and ~ \mu \in \La_{n}\backslash \La_{n-4}$~we have$$\Ext_{\scriptsize{A_{n}}}^{1}(\De_{n}(\la),\De_{n}(\mu) )=0,$$
Then each of the endomorphism algebras~$A_{n}$ in~$\A$~is semi-simple.
\end{thm}

To prove Theorem \ref{semis}, we first need the following lemmas.

\begin{lemma}\label{hufan}
Let $A$ be a $K$-algebra. Suppose a~$K$-algebra~$B$~is a semi-simple subalgebra of~$A$. Then for any~$A$-module~$M$~and~$B$-module~$N$, we have
\[
\Ext_{\scriptsize{A}}^{i}(\ind (N),M)\left\{
\begin{array}{ll}
\simeq  \Hom_{\scriptsize{B}}(N, \res(M)) &\mbox{if~$i=0$,}\\
=0&\mbox{if~$i>0$.}
\end{array}
\right.
\]

\end{lemma}
\proof
Let~$M\rightarrow Q^{\bullet} $~be an injective resolution of~$M$~in~$A\mod$. Applying~$\res$, due to~$B$~is semi-simple, it follows that~$\res M\rightarrow (\res (Q))^{\bullet} $ is an injective resolution of ~$\res (M)$~in~$B \mod$
and
\begin{eqnarray*}
\Ext_{\scriptsize{B}}^{i}(N, \res (M))& = &H^{i}(\Hom_{\scriptsize{B}}(N, (\res (Q))^{\bullet}))\\
&= & H^{i}(\Hom_{\tiny{B}}(N, \Hom(A_{\scriptsize B},Q^{\bullet} )))\\
&\simeq & H^{i}(\Hom_{\scriptsize{A}}(A\o_{\tiny{B}} N, Q^{\bullet} ))\\
&= & \Ext_{\scriptsize{A}}^{i}(\ind (N), M).
\end{eqnarray*}

Since~$B$~is semi-simple, we have~$\Ext_{\scriptsize{B}}^{i}(N, \res (M))=0$ for~$i>0$.
 The lemma follows.
$\qed$

Martin~and~Woodcock~\cite{MW}~studied standard modules of quasi-hereditary algebras in terms of induction and restriction, the following result is a generalizition of ~\cite[Proposition 3.5]{MW}.
\begin{lemma}\label{ext}
Let~$\A$~be a hereditary cellular category. Suppose~$m, n\in \mathbb{N}$~satisfying~$\Hom_{\A}(m,n)\neq 0$. Then for each~$\la, \mu \in \La_{(n,m)}$~and~$i\geq 0$, we have
$$\Ext_{\scriptsize{A_{m}}}^{i}(\De_{m}(\la), \De_{m}(\mu))\simeq \Ext_{\scriptsize{A_{n}}}^{i}(\De_{n}(\la), \De_{n}(\mu)).$$

\end{lemma}
\proof Without loss of generality, we may suppose~$n\leq m$.
Let~$P^{\bullet}\rightarrow \De_{n}(\la)$~be a projective resolution of~$\De_{n}(\la)$~in~$A_{n}\mod$. Denoted by~$F:=\Hom_{\A}(m,n)\o_{\tiny{A_{n}}}-$, due to~Proposition \ref{hanzi2},~$F$~is an idempotent embedding functor.
Because by assumption
~$\La_{n}$~is a saturated subset of~$\La_{m}$, due to~\cite[Proposition A 3.2]{Do}, for any~$i>0$, we have~$L_{i}F\De_{n}(\la)=0$, where~$L_{i}F$~denotes the $i$th left derived functor. Moreover, due to Lemma~\ref{hanzi1}, it follows that
~$(\Hom_{\A}(m,n)\o P^{\bullet})\rightarrow \De_{m}(\la)$~is a projective resolution of~$\De_{m}(\la)$~in~$A_{m}\mod$. 
Let~$e$~be an idempotent  of~$A_{m}$~satisfying
~$_{\scriptsize{A_{m}}}\Hom_{\A}(m,n)\simeq(A_{m}) e$,
we have
\begin{eqnarray*}
\Ext_{\scriptsize{A_{m}}}^{i}(\De_{m}(\la), \De_{m}(\mu))& = &H^{i}(\Hom_{\scriptsize{A_{m}}}((\Hom_{\A}(m,n)\o P^{\bullet}), \De_{m}(\mu)))\\
&\simeq & H^{i}(\Hom_{\scriptsize{A_{m}}}(((A_{m}) e\o P^{\bullet}), \De_{m}(\mu)))\\
&\simeq & H^{i}(\Hom_{\scriptsize{A_{n}}}( P^{\bullet}, e(\De_{m}(\mu))))\\
&= & \Ext_{\scriptsize{A_{n}}}^{i}(\De_{n}(\la), \De_{n}(\mu)).
\end{eqnarray*}
$\qed$

We are now in the position to give a proof of Theorem~\ref{semis}.

\proof
By the assumption that~$A_{n}$~is quasi-hereditary and due to Theorem~\ref{cao}, we need to prove,
for all~$\la, \mu \in \La_{n}$~with~$\mu <\la$, we have
$$\Ext_{\scriptsize{A_{n}}}^{1}(\De_{n}(\la),\De_{n}(\mu) )=0.$$
Thus we always assume that~$\mu <\la$.

We use induction on~$n$. Assume that~$\la \in \La_{n-2}$, then it is clear that~$\mu \in \La_{n-2}$. Due to Lemma~\ref{ext},
$$\Ext_{\scriptsize{A_{n}}}^{1}(\De_{n}(\la), \De_{n}(\mu))\simeq \Ext_{\scriptsize{ A_{n-2}}}^{1}(\De_{n-2}(\la), \De_{n-2}(\mu)).$$
According to the induction hypothesis, the right  side of the above  vanishes. Hence so is the left.

Assume that~$\la \in \La_{n}\backslash \La_{n-2}$~, $ \mu \in \La_{n}\backslash \La_{n-4}$,  it follows directly from the assumptions  that
$$\Ext_{\scriptsize{A_{n}}}^{1}(\De_{n}(\la),\De_{n}(\mu) )=0.$$

Assume now that~$\la \in \La_{n}\backslash \La_{n-2}$, $ \mu \in \La_{n-4}$,
by axiom~$\mbox{\rm (}A4\mbox{\rm )}$, there exists~$\tau \in \La_{n-1}\backslash \La_{n-3}$~
such that there is an exact sequence
\begin{equation}\label{K}
0\rightarrow K \rightarrow\ind_{n-1}(\De_{n-1}(\tau))\rightarrow \De_{n}(\la)\rightarrow 0.
\end{equation}
and
$$\supp_{n}(\ind_{n-1}(\De_{n-1}(\tau)))\subseteq (\La_{n-2}\backslash \La_{n-4} )\cup (\La_{n}\backslash \La_{n-2} ).$$
By~(\ref{K}), we have the exact sequence
$$
0\rightarrow \Hom_{\scriptsize{A_{n}}}(\De_{n}(\la),\De_{n}(\mu))
$$
\begin{equation}\label{long}
\rightarrow\Hom_{\scriptsize{A_{n}}}(\ind_{n-1}(\De_{n-1}(\tau)),\De_{n}(\mu) )\rightarrow\Hom_{\scriptsize{A_{n}}}(K, \De_{n}(\mu))\end{equation}
$$
\rightarrow\Ext_{\scriptsize{A_{n}}}^{1}(\De_{n}(\la),\De_{n}(\mu) )
\rightarrow \Ext_{\scriptsize{A_{n}}}^{1}(\ind_{n-1}(\De_{n-1}(\tau)),\De_{n}(\mu) ).$$

Due to~$\mbox{\rm (}A3\mbox{\rm )}$,  it follows that~$\supp_{n}(\res_{n}(\De_{n}(\mu)))\subseteq \La_{n-3}$.

Furthermore, by the induction hypothesis,~$ A_{n-1}$~is semi-simple. Hence~$\res_{n}(\De_{n}(\mu))\simeq \bigoplus\limits_{i}\De_{n-1}(\upsilon_{i})$~with~some~$\upsilon_{i}\in \La_{n-3}$.
According to Lemma~\ref{hufan}, we also have $$\Ext_{\scriptsize{A_{n}}}^{1}(\ind_{n-1}(\De_{n-1}(\tau)),\De_{n}(\mu))=0.$$
and
\begin{eqnarray*}
&&\Hom_{\scriptsize{A_{n}}}(\ind_{n-1}(\De_{n-1}(\tau)),\De_{n}(\mu) )\\
& = &\Hom_{\scriptsize{ A_{n-1}}}(\De_{n-1}(\tau),\res_{n}(\De_{n}(\mu)) )\\
&= & \Hom_{\scriptsize{ A_{n-1}}}(\De_{n-1}(\tau),\bigoplus\limits_{i}\De_{n-1}(\upsilon_{i}) )\\
&= & 0.
\end{eqnarray*}
Consequently, in the long sequence (\ref{long}), for ~$\la \in \La_{n}\backslash \La_{n-2}$, $ \mu \in \La_{n-4}$, $ \Hom_{\scriptsize{A_{n}}}(\De_{n}(\la),\De_{n}(\mu))=0$.

To prove ~$\Ext_{\scriptsize{A_{n}}}^{1}(\De_{n}(\la),\De_{n}(\mu) )=0$, it is sufficient to show  $\Hom_{\scriptsize{A_{n}}}(K,\De_{n}(\mu))=0$ in the long sequence (\ref{long}).

Firstly,  for~$\la \in \La_{n-2}\backslash\La_{n-4}$, $\mu \in \La_{n-4}$, it is clear that
$$~ \Hom_{\scriptsize{A_{n}}}(\De_{n}(\la),\De_{n}(\mu))=\Hom_{\scriptsize{ A_{n-2}}}(\De_{n-2}(\la),\De_{n-2}(\mu))=0.$$
Thus, for all~$\la \in \La_{n}\backslash\La_{n-4}$, $\mu \in \La_{n-4}$ we have~$\Hom_{\scriptsize{A_{n}}}(\De_{n}(\la),\De_{n}(\mu))=0$.

Moreover, because
$\supp_{n}(K)\subseteq  \La_{n}\backslash\La_{n-4}$, this implies~$\Hom_{\scriptsize{A_{n}}}(K,\De_{n}(\mu))=0$.
Hence~$\Ext_{\scriptsize{A_{n}}}^{1}(\De_{n}(\la),\De_{n}(\mu) )=0$.
This finishes the proof..
$\qed$

There are a large number of concrete algebras in our axiom scheme. We briefly list some as follows.

\textbf{Example 1}~~$\mbox{\rm (}1\mbox{\rm )}$~~ The Temperley-Lieb categoy~$\mathcal{TL}(\de)$~with~$\de$ not a root of unity. See Section 4.2 for the definition. For $n\in \mathbb{N}$, the poset is~$\La_{n}=\{0~or~1,\ldots, n-4, n-2, n \}$. For each~$0\leq i<n$, there is  a short exact sequence
$$0\rightarrow \De_{n+1}(i-1) \rightarrow\ind_{n}(\De_{n}(i))\rightarrow \De_{n+1}(i+1)\rightarrow 0.$$

For more details on restriction and induction for Temperley-Lieb algebras, the reader is referred to ~\mbox{\rm{\cite{Wes2,RSA}}}.

~$\mbox{\rm (}2\mbox{\rm )}$ The Brauer diagram categoy~$\mathcal{B}(\de)$~with~$\de\neq 0$. For the definition we refer the reader to \cite{LZ} for details.

 For $n\in \mathbb{N}$, the endomorphism algebra of $n$ is the~Brauer~algebra~$B_{n}(\de)$, due to~Cox~\mbox{\rm{\cite{CVP}}}, it is
quasi-hereditary whenever it is characteristic zero or characteristic zero $p>n$.

By~\mbox{\rm{\cite[Theorem~4.1~and~Corollary~6.4]{DWH}}}~and~\mbox{\rm{\cite[Proposition~2.7]{CVP}}}, in arbitrary  characteristic, for each~$\la \in \La_{n}$,  there is  a short exact sequence
$$0\rightarrow \bigoplus\limits_{\tau \lhd \la}\De_{n+1}(\tau) \rightarrow\ind_{n}(\De_{n}(\la))\rightarrow \bigoplus\limits_{\mu \rhd \la}\De_{n+1}(\mu)\rightarrow 0,$$
where~$\La_{n}$~is the indexing set of~$B_{n}(\de)$,  $\tau \lhd \la$~represents that  the young tableau $\la$  is obtained by adding a box from  the~young tableau~$\tau$ and  $\mu \rhd \la$~represents that  the young tableau $\la$ is obtained by removing a box from  the~young tableau~$\mu$.

For similar resaults on~partition algebras~and~blob algebras and their categorical analogues, we refer the reader to~\mbox{\rm{\cite{ CGM,Martin, MR, MW}}}.

\section{Algebraic structures on Grothendieck groups of tower categories}
A number of diagram categories of great interest, such as Temperley-Lieb categories, Brauer diagram categories and so on,
are cellular. The common feature of these categories is that they admit a natual tensor structure, the `juxtaposition' of
two diagrams, that is, the tensor product $A\o B$ of two diagrams is obtained by putting the diagram of $A$ on the left of the diagram of $B$.
Another feature is that  $ A_{0} \cong K \cong A_{1}$. So we consider the following axiom in our framework and  equip the hereditary cellular category with a
 tensor product.  

$\mbox{\rm (}A1^\prime \mbox{\rm )}$ {\bf $ A_{0} \cong K$, $A_{1} \cong K$ and there is an external multiplication $\rho_{m,n}: A_{m} ~\otimes
A_{n}  \rightarrow  A_{m+n},$~for all $m,n\geqslant 0,$ such that}
 
 ~~~~~~{\bf (a)~For all $m$ and $n$,  $\rho_{m,n}$ is an injective homomorphism of algebras, sending $1_{m}\otimes 1_{n} $~to
$1_{m+n} $};

~~~~~~{\bf (b)~ $\rho$ is associative, that is, $\rho_{l+m,n}~\cdot~(\rho_{l,m}\otimes 1_{n})=\rho_{l,m+n}~\cdot~(1_{n}\otimes \rho_{m,n} )$, for all $l,m,n$.}

Obviously, $\mbox{\rm (}A1\mbox{\rm )}$, in Definition \ref{deftower}, is the special case 
for $m=1$ in $\mbox{\rm (}A1^\prime \mbox{\rm )}$. The following definition is obtained by replacing $\mbox{\rm (}A1\mbox{\rm )}$ with $\mbox{\rm (}A1^\prime\mbox{\rm )}$.

\begin{defin}\label{defTTtower}
A  hereditary cellular category is called a  tensor cellular tower category (TCTC) if it satisfies
~$\mbox{\rm (}A1^\prime\mbox{\rm )}$ and~$\mbox{\rm (}A2\mbox{\rm )}$--~$\mbox{\rm (}A4\mbox{\rm )}$.
\end{defin}

Let $\A$  be a TCTC.  Denote by $G_{0}(A_{n})$ the Grothendieck group of the endomorphism algebra of object $n$.
Let $G_{0}(\A):=\oplus_{n\geqslant 0}G_{0}(A_{n})$. The second reason allowing tensor product
is that one can consider the
algebraic structures on Grothendieck groups. 
In \cite{BL}, Bergeron and Li introduced a general notation of a tower of algebras by way of axioms,
which guarantee that the Grothendieck groups of 
a tower of algebras $\oplus_{n\geqslant 0}G_{0}(A_{n})$
can be a pair of graded dual Hopf algebras.
The tower of symmetric group algebras and their quantum deformations (such as Hecke algebras and Hecke-Clifford algebras and so on) are typical examples in their notation.  The common feature of those examples is that they admit Mackey's formula, which is
just the compatibility relation between the multiplication and comultiplication on Grothendieck groups. Then they gave an analogue of Mackey's formula (see  equality  (5) below) as an axiom in their framework.
We shall see that the analogue of Mackey's formula never holds in  the quasi-hereditary case. The Grothendieck groups, however, have algebra and coalgebra structures under  certain conditions.

In this section, the Grothendieck groups for TCTC
 is studied
in 4.1.
As a typical example of TCTC, we study  Temperley-Lieb categories and obtain  the structure constants of the multiplication on the Grothendieck groups in 4.2.

\subsection{Induction and restriction functors on $G_{0}(\A)$}
Let $\A$  be a TCTC, we next consider the  induction and restriction on $G_{0}(\A)$ in terms of tensor products.

For $m,n \in \mathbb{N}$, let $M$ be a left $A_{ m}$-module, and  $N$ be a left $A_{ n}$-module. Recall that the tensor product $M\otimes_{K} N$ is a left
$A_{ m} \otimes_{K} A_{ n}$-module with the action $(a \otimes b)\cdot (w\otimes u)=aw \otimes bu$ for $a\in A_{ m}, b\in A_{ n}, w\in M$ and $u\in N$.

Define the inductions on $G_{0}(\A)$ as follows:

$$i_{m,n} :   G_{0}(A_{m})\otimes_{\mathbb{Z}} G_{0}(A_{n})\rightarrow G_{0}(A_{m+n})$$
$$[M]\otimes [N] \mapsto [\ind_{\scriptsize{A_{ m} \otimes A_{ n}}}^{\scriptsize{A_{ m+n}}}( M\otimes N)],$$  where
$$\ind_{\scriptsize{ A_{m} \otimes  A_{n}}}^{\scriptsize{ A_{m+n}}} (M\otimes N)= A_{m+n} \bigotimes_{ A_{m} \otimes  A_{n}}(M\otimes N)$$ 
$$=\frac{ A_{m+n} \otimes(M\otimes N)}{<a\otimes [(b\otimes c)(w\otimes u)]-[a\rho_{m,n}(b\otimes c)]\otimes w\otimes u>},
$$ 
for $a\in  A_{m+n}, ~ b\in  A_{m},~ c\in  A_{n},~w\in M$~and~$u\in N$.

Also define
$$r_{k,l}:  G_{0}( A_{n})\rightarrow G_{0}( A_{k})\bigotimes_{\mathbb{Z}}G_{0}( A_{l})~with~k+l=n$$
$$ [N] \mapsto [\res_{\scriptsize{ A_{k} \otimes  A_{l}}}^{\scriptsize{ A_{n}}}~(N)],$$ 
where $\res_{\scriptsize{ A_{k} \otimes  A_{l}}}^{\scriptsize{ A_{n}}}~(N)~=~\Hom_{\scriptsize{ A_{n}}}( A_{n},N)$ is an 
$ A_{k} \otimes  A_{l}$-module with the action defined by $((b\otimes c)\cdot f)(a)=f(a\rho_{k,l}(b\otimes c))$ for 
$a\in A_{n},~ b\in  A_{k},~ c\in  A_{l}$~and~$f\in \Hom_{\scriptsize{ A_{n}}}( A_{n},N).$

In \cite{BL}, the following condition is as an axiom  to insure that the  maps $i$~and~$r$~are well defined on for a tower of algebras:

$ A_{m+n}$ is a two-sided projective $A_{m} \otimes A_{n}$-module with the action defined by $a\cdot(b\otimes c)=a\rho_{m,n}(b\otimes c)$~and~$(b\otimes c)\cdot a=\rho_{m,n}(b\otimes c)a$
for all
$m,n \geqslant 0, a\in  A_{m+n},~ b\in A_{m},~ c\in A_{n}.$

Thus, we immediately have

\begin{prop}
Let $\A$ be a TCTC. Suppose $ A_{m+n}$ is a two-sided projective $A_{m} \otimes A_{n}$-module for any $m,n \in \mathbb{N}$. Then the maps $i$~and~$r$~are well defined on $G_{0}(\A)$.
\end{prop}

Further, we construct the multiplication and comultiplication by $i$ and $r$ and define the unit and counit on $G_{0}(\A)$~as follows:

$$\pi~: G_{0}(\A)\otimes_{\mathbb{Z}}G_{0}(\A)\rightarrow G_{0}(\A)$$
$$~~~~~~~~~~~~~~~\textrm{where}~\pi|_{G_{0}(\scriptsize{ A_{k}})\otimes G_{0}(\scriptsize{ A_{l}} )}=i_{k,l}$$
$$\Delta~: G_{0}(\A) \rightarrow G_{0}(\A)\otimes_{\mathbb{Z}}G_{0}(\A) $$
$$~~~~~~~~~~\textrm{where}~\Delta|_{G_{0}(\scriptsize{A_{n}} )}=\sum_{k+l=n}r_{k,l}$$
$$\mu~: \mathbb{Z}\rightarrow G_{0}(\A)~~~~~~~~~~~~~~~~~~~~~$$
$$~~~~~~~~~~~~~~~~~~~~~~~~~~\textrm{where}~\mu(a)=a[K]\in {G_{0}( A_{0})}, \textrm{for}~ a\in \mathbb{Z}$$
$$\epsilon~:G_{0}(\A)\rightarrow  \mathbb{Z}~~~~~~~~~~~~~~~~~~~~~$$
\[~~~~~~~~~~~~~~~~~~~~~~~~~~~~~~
\textrm{where}~\epsilon([M])=\left\{
\begin{array}{ll}
a &\mbox{if~$[M]=a[K]\in G_{0}( A_{0})$,~where~$a\in\mathbb{Z}$,}\\
0&\mbox{otherwise.}
\end{array}
\right.
\]

Also, Bergeron and Li \cite[Theorem 3.5]{BL} proved that the associativity of $\pi$, the unitary property of $\mu$, the coassociativity of $\Delta$~and the counitary property of $\epsilon$, which imply that $ (G_{0}(\A), \pi, \mu)$ is an algebra and $(G_{0}(\A), \Delta, \epsilon)$ is a coalgebra.
Thus, we have

\begin{thm}\label{alc}
Let $\A$ be a TCTC. Suppose $ A_{m+n}$ is a two-sided projective $A_{m} \otimes A_{n}$-module for any $m,n  \in \mathbb{N}$. Then $ (G_{0}(\A), \pi, \mu)$~is an algebra and  $(G_{0}(\A), \Delta, \epsilon)$ is a coalgebra.
\end{thm}

As a special case in Theorem \ref{alc}, it is easy to see that if each $A_{n}$ is semi-simple then $ (G_{0}(\A), \pi, \mu)$~is an algebra and  $(G_{0}(\A), \Delta, \epsilon)$ is a coalgebra.

It is natural to ask that whether $G_{0}(\A)$ is a bialgebra? This means that whether the equality $$\Delta (\pi([M]\otimes [N]))=\pi (\Delta ([M])\otimes \Delta ([N])) \eqno(4)$$ holds, that is, $\Delta$ is an algebra homomorphism.
 
To check this, however, it  first needs to define a
reasonable  multipulication  on $\Delta ([M])\otimes \Delta ([N])$. To do this, note that there is a natural manner in terms of the twisted tensor product which  is called the \textit{twisted induction} in \cite{BL}.
 
Let $k=t+s$, define
$$ G_{0}( A_{t}\otimes  A_{m-t})\bo_{\mathbb{Z}} G_{0}( A_{s}\otimes  A_{n-s})\rightarrow G_{0}( A_{k}\otimes  A_{m+n-k})$$
$$[M_{1}\o M_{2}]\otimes [N_{1}\o N_{2}] \mapsto [\widetilde{\ind}_{\scriptsize{ A_{t} \o  A_{m-t}\o A_{s} \o  A_{n-s}}}^{\scriptsize{ A_{k}\otimes  A_{m+n-k}}} ((M_{1}\o M_{2})\otimes (N_{1}\o N_{2}))],$$  

where
$$\widetilde{\ind}_{\scriptsize{ A_{t} \o  A_{m-t}\o A_{s} \o  A_{n-s}}}^{\scriptsize{ A_{k}\otimes  A_{m+n-k}}} ((M_{1}\o M_{2})\otimes (N_{1}\o N_{2}))$$
$$=( A_{k}\otimes  A_{m+n-k})\widetilde{\bo}_{\scriptsize{ A_{t} \o  A_{m-t}\o A_{s} \o  A_{n-s}}}((M_{1}\o M_{2})\otimes (N_{1}\o N_{2}))$$
This means
$$(a\o b)\o [(c_{1}\o c_{2})\cdot(w_{1}\cdot w_{2})\o (d_{1}\o d_{2})\cdot(u_{1}\cdot u_{2})]$$
$$\equiv [a\rho_{t,s}(c_{1}\o d_{1})\o b\rho_{m-t, n-s}(c_{2}\o d_{2})]\o (w_{1}\o u_{1} \o w_{2} \o u_{2}).$$

Therefore equality (4) can be interpreted as the following equality:
$$[\res_{\scriptsize{ A_{k}\otimes  A_{m+n-k}}}^{\scriptsize{ A_{m+n}}}\ind _{\scriptsize{A_{n}\otimes A_{m}}}^{\scriptsize{ A_{m+n}}}(M\o N)]$$
$$=\sum_{t+s=k}[\widetilde{\ind}_{\scriptsize{ A_{t} \o  A_{m-t}\o A_{s} \o  A_{n-s}}}^{\scriptsize{ A_{k}\otimes  A_{m+n-k}}}(\res_{\scriptsize{ A_{t}\o  A_{m-t}}}^{\scriptsize{A_{m}}}(M) \o \res_{\scriptsize{ A_{s}\o  A_{n-s}}}^{\scriptsize{A_{n}}}(N))]\eqno(5)$$
for all $0<k<m+n$,  an $A_{m}$-module $M$ and  an $A_{n}$-module $N$.

In the rest of the subsection, we  verify that equality (5) is not true for Temperley-Lieb categories.  Let~$\mathcal{TL}(\delta)$~be a $\mathcal{TL}$-category with $\de$ not a root of unity,  and let $\De_{n}(r)$ be the cell module  of the endomorphism algebra $TL_{n}$ corresponding to the index $r$.

It is known that (see \cite[ Corollary 4.6 and Corollary 6.6]{RSA})
$$\res_{A_{n-1}}^{A_{n}}(\De_{n}(p))\cong \De_{n-1}(p)\oplus \De_{n-1}(p-1) ~~~(as~a~ TL_{n-1}-module) $$ and
$$\ind_{A_{n}}^{A_{n+1}}(\De_{n}(p))\cong \De_{n+1}(p+1)\oplus \De_{n+1}(p) ~~~(as~a~ TL_{n+1}-module) .$$

Then, for cell modules $\De_{n}(p)$~and~$ \De_{1}(0)$, the left side of (5) equals:
\begin{eqnarray*}  
&&[\res_{\scriptsize{A_{n}\otimes A_{1}}}^{\scriptsize{ A_{n+1}}}\ind _{\scriptsize{A_{n}\otimes A_{1}}}^{\scriptsize{ A_{n+1}}}(\De_{n}(p)\o \De_{1}(0))]\\
& = &[\res_{\scriptsize{A_{n}\otimes A_{1}}}^{\scriptsize{ A_{n+1}}}(\De_{n+1}(p+1)\o \De_{1}(0))\oplus(\De_{n+1}(p)\o \De_{1}(0))]\\
&= &[\res_{\scriptsize{A_{n}\otimes A_{1}}}^{\scriptsize{ A_{n+1}}}(\De_{n+1}(p+1)\o \De_{1}(0))]
+[ \res_{\scriptsize{A_{n}\otimes A_{1}}}^{\scriptsize{ A_{n+1}}}(\De_{n+1}(p)\o \De_{1}(0))]\\
&=& [\res_{\scriptsize{A_{n}\otimes A_{1}}}^{\scriptsize{ A_{n+1}}}(\De_{n+1}(p+1))]
+[\res_{\scriptsize{A_{n}\otimes A_{1}}}^{\scriptsize{ A_{n+1}}}(\De_{n+1}(p))]\\
&=& [\De_{n}(p+1)] + 2[\De_{n}(p)]  +[\De_{n}(p-1)] .
\end{eqnarray*}
On the other hand, the right side of (5) is:
\begin{eqnarray*}
&&\sum_{t+s=1}[\widetilde{\ind}_{\scriptsize{ A_{t} \o  A_{n-t}\o A_{s} \o  A_{1-s}}}^{\scriptsize{A_{1}\otimes A_{n}}}(\res_{\scriptsize{ A_{t}\o A_{n-t}}}^{\scriptsize{A_{n}}}(\De_{n}(p)) \o \res_{\scriptsize{ A_{s}\o  A_{1-s}}}^{\scriptsize{A_{1}}}(\De_{1}(0)) )]\\
& = &[ \widetilde{\ind}_{\scriptsize{ A_{0} \o A_{n}\o A_{1} \o  A_{0}}}^{\scriptsize{A_{1}\otimes A_{n}}}(\res_{\scriptsize{ A_{0}\o A_{n}}}^{\scriptsize{A_{n}}}(\De_{n}(p)) \o \res_{\scriptsize{A_{1}\o  A_{0}}}^{\scriptsize{A_{1}}}(\De_{1}(0)) )]\\
&&+[\widetilde{\ind}_{\scriptsize{A_{1} \o  A_{n-1}\o A_{0} \o A_{1}}}^{\scriptsize{A_{1}\otimes A_{n}}}(\res_{\scriptsize{A_{1}\o  A_{n-1}}}^{\scriptsize{A_{n}}}(\De_{n}(p)) \o \res_{\scriptsize{ A_{0}\o  A_{1}}}^{\scriptsize{A_{1}}}(\De_{1}(0)) )]\\
&=& [\ind_{\scriptsize{A_{n}}}^{\scriptsize{A_{n}}}(\De_{n}(p)\o \De_{1}(0))]+ [\ind_{\scriptsize{ A_{n-1}\o A_{1}}}^{\scriptsize{A_{n}}}(\De_{n-1}(p)\oplus \De_{n-1}(p-1))] \\
&=&[ \ind_{\scriptsize{A_{n}}}^{\scriptsize{A_{n}}}(\De_{n}(p)\o \De_{1}(0))]+[ \ind_{\scriptsize{ A_{n-1}}}^{\scriptsize{A_{n}}}(\De_{n-1}(p)) \oplus\ind_{\scriptsize{ A_{n-1}}}^{\scriptsize{A_{n}}} (\De_{n-1}(p-1))]\\
&=& [\De_{n}(p)] + [\De_{n}(p+1) \oplus \De_{n}(p) \oplus \De_{n}(p) \oplus \De_{n}(p-1) ]\\
&=&  [\De_{n}(p+1)] + 3[\De_{n}(p)]  +[\De_{n}(p-1)] .
\end{eqnarray*}
Consequently, the equality (5) do not hold.
Hence  $ (G_{0}(\mathcal{TL}(\delta)), \pi, \mu, \Delta, \epsilon)$~is not a bialgebra.

\subsection{Algebraic structures on the Grothendieck groups of $\mathcal{TL}$-categories}
In this section, we shall study the multiplication in the Grothendieck groups of semi-simple $\mathcal{TL}$-categories, and have the following result.

\begin{thm}\label{cons}
Let~$\mathcal{TL}(\delta)$~be a $\mathcal{TL}$-category with $\de$ not a root of unity. Then: 

For each $m,n \in \mathbb{N}$, 
$$i_{m,n}[[\Delta_{m}(p)]\otimes [\De_{n}(q)]] =\sum_{0 \leqslant r \leqslant [(m+n)/2]}a_{(p, q, r)}^{(m|n)}[\Delta_{m+n}(r)],$$  
where
\begin{eqnarray*}
a_{(p, q, r)}^{(m|n)}&=&\left\{
\begin{array}{ll}
1 &\mbox{if~$p+q\leqslant r$~ and ~$m-s\geqslant 2p$~and ~$m-s\geqslant 2q$,}\\
0 &\mbox{others.}
\end{array}
\right. ,
\end{eqnarray*}
with $s=r-(p+q)$ and $\Delta_{k}(l)$, $0\leqslant l \leqslant [k/2]$, the cell module of $TL_{k}$. 
\end{thm}

It is known that Temperley-Lieb algebra $TL_{n}$ is semi-simple if $\delta$ is not a root of unity, and hence all cell modules of $TL_{n}$ form a complete set of non-isomorphism simple $TL_{n}$-modules by Theorem \ref{sem}. Moreover, all $\Delta_{m}(p)\otimes_{K} \Delta_{n}(q)$ for $0\leqslant p \leqslant [m/2]$ and $0\leqslant q \leqslant [n/2]$ 
form a complete set of non-isomorphism simple modules
for semi-simple algebra $TL_{m}\o_{K} TL_{n}$.

\begin{figure}[h]
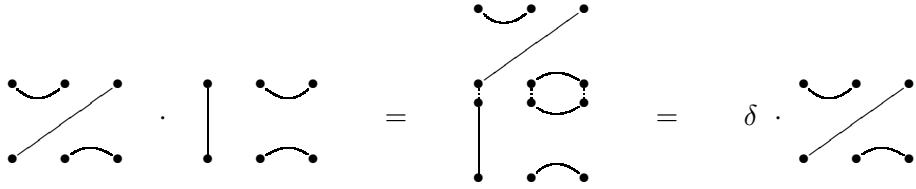

\[
\xy
(-8,-5.5)*{\scriptstyle\bullet}="1s"; 
(-1,-5.5)*{\scriptstyle\bullet}="2s"; 
(6,-5.5)*{\scriptstyle\bullet}="3s";
(-8,4.5)*{\scriptstyle\bullet}="1d"; 
(-1,4.5)*{\scriptstyle\bullet}="2d"; 
(6,4.5)*{\scriptstyle\bullet}="3d";
"1d"; "2d" **\crv{(-6,2) & (-3.5,2)};
"3d"; "1s" **\dir{-};
"2s"; "3s" **\crv{(1,-3.5) & (3.5,-3.5)};
\endxy
~~~\cdot~~~
\xy
(13,-5.5)*{\scriptstyle\bullet}="4s";
(20,-5.5)*{\scriptstyle\bullet} ="5s";
(27,-5.5)*{\scriptstyle\bullet}="6s";
(13,4.5)*{\scriptstyle\bullet}="4d";
(20,4.5)*{\scriptstyle\bullet} ="5d";
(27,4.5)*{\scriptstyle\bullet}="6d";
"5d"; "6d" **\crv{(22.5,2) & (25,2)};
"4d"; "4s" **\dir{-};
"5s"; "6s" **\crv{(22.5,-3.5) & (25,-3.5)};
\endxy
~~~~~~
=
~~~~~~
\xy
(-8,2)*{\scriptstyle\bullet}="1"; 
(-1,2)*{\scriptstyle\bullet}="2"; 
(6,2)*{\scriptstyle\bullet}="3";
(6,-8)*{\scriptstyle\bullet}="12"; 
(-1,-8)*{\scriptstyle\bullet}="13"; 
(-8,-8)*{\scriptstyle\bullet}="10"; 
(-8,4.5)*{\scriptstyle\bullet}="1s"; 
(-1,4.5)*{\scriptstyle\bullet}="2s"; 
(6,4.5)*{\scriptstyle\bullet}="3s";
(-8,14.5)*{\scriptstyle\bullet}="1d"; 
(-1,14.5)*{\scriptstyle\bullet}="2d"; 
(6,14.5)*{\scriptstyle\bullet}="3d";
"1"; "10" **\dir{-};
"1d"; "2d" **\crv{(-6,12) & (-3.5,12)};
"13"; "12" **\crv{(1,-5.5) & (3.5,-5.5)};
"3d"; "1s" **\dir{-};
"2s"; "3s" **\crv{(1,6.5) & (3.5,6.5)};
"2"; "3" **\crv{(1,0) & (3.5,0)};
"1"; "1s" **\dir{.};
"2"; "2s" **\dir{.};
"3"; "3s" **\dir{.};

\endxy
~~~~~~
=
~~~~~~
\delta~\cdot~
\xy
(-8,-5.5)*{\scriptstyle\bullet}="1s"; 
(-1,-5.5)*{\scriptstyle\bullet}="2s"; 
(6,-5.5)*{\scriptstyle\bullet}="3s";
(-8,4.5)*{\scriptstyle\bullet}="1d"; 
(-1,4.5)*{\scriptstyle\bullet}="2d"; 
(6,4.5)*{\scriptstyle\bullet}="3d";
"1d"; "2d" **\crv{(-6,2) & (-3.5,2)};
"3d"; "1s" **\dir{-};
"2s"; "3s" **\crv{(1,-3.5) & (3.5,-3.5)};
\endxy
\]
	\caption{the multiplication of $3$-diagrams}
	\label{fig1}
\end{figure}

Recall that an $n$-diagram consist of two rows of $n$ dots in which each dot is joined to just one other dot and none of the joins intersect when drawn in the rectangle defined by the two rows of 
$n$ dots. All such $n$-diagrams form a cellular basis of $TL_{n}$ and the multiplication of the basis is given by concatenating two $n$-diagrams, that is, stacking the first diagram  on top of the second diagram, matching the relevant bottom and top vertices, and replacing all closed cycles by a factor
$\de^{r}$ with $r$ the number of closed cycles.   See Figure \ref{fig1} for an example for $n=3$.

Due to \cite{GL}, the cell module $\De_{n}(r)$ of $TL_{n}$ is spanned by all $(n, r)$-cap diagrams, which is the `half-diagrams' obtained from an $n$-diagram by cutting horizontally down the middle. Figure \ref{fig2}(2) illustrates
an (11,4)-cap diagram in which the arcs $\{3,8\},\{4,5\}$, $\{6,7\}$ and $\{10,11\}$ are caps and vertices $1,2,9$ are called single points.

For further informations on the diagrammatic and algebraic defintions of Temperley-Lieb algebras and their  cell modules, we refer the reader to \cite{GL,RSA,Wes2}.

To prove Theorom \ref{cons}, we first regard a $(m+n)$-diagram (and a cap diagram) as a `walled' diagram  in Definition \ref{wall}.
By using such `walled' diagrams, we then give a composition series of $_{TL_{m}\o TL_{n}}\De_{m+n}(r)$ in Proposition \ref{keyp}. Finally, Theorom \ref{cons} holds true directly  from  the following equalities.
$$i_{m,n}[[\Delta_{m}(p)]\otimes [\De_{n}(q)]]=  [\ind_{\scriptsize{A_{ m} \otimes A_{ n}}}^{\scriptsize{A_{m+n}}}[ \Delta_{m}(p)\otimes \Delta_{n}(q)]]=\sum_{r}a_{(p,q,r)}^{(m|n)}[\Delta_{m+n}(r)],$$  
and
\begin{align}
	a_{(p,q,r)}^{(m|n)}&=&\dim \Hom[\ind_{\scriptsize{A_{ m} \otimes A_{ n}}}^{\scriptsize{A_{m+n}}} [\Delta_{m}(p)\otimes \Delta_{n}(q)],   \Delta_{m+n}(r)] \notag\\
& = &\dim \Hom[_{\scriptsize{A_{ m} \otimes A_{ n}}} \Delta_{m}(p)\otimes \Delta_{n}(q),  _{\scriptsize{A_{ m} \otimes A_{ n}}}  \Delta_{m+n}(r)]. \tag{6}
\end{align}

The following observations is the key  to our proof.
Indeed, the juxtaposition 
(see the beginning of the section 4) of an $m$-diagram and an $n$-diagram can be  viewed as a `walled' diagram if we imagine that there is a wall between the two diagrams, and we herein call such diagram an $(m|n)$-diagram. It is clear that all such diagrams form a basis of $TL_{m}\o TL_{n}$. Figure \ref{fig2}(1) demonstrates the juxtaposition 
of an $6$-diagram in $TL_{6}$ and an $5$-diagram in $TL_{5}$.

Similarly, an~$(m+n, r)$-cap diagram in $_{TL_{m}\o TL_{n}}\De_{m+n}(r)$ can be  viewed as a `walled' cap diagram if we imagine that there is a wall which separates the vertices $m, m+1$. Figure \ref{fig2}(2) illustrates that
an (11,4)-cap diagram can be viewed as a `walled' diagram
when $\De_{11}(4)$ is restricted as a $TL_{6}\o TL_{5}$- module.

Due to the above observations, we define

\begin{defin}\label{wall}
Let~$\mathcal{TL}(\delta)$~be a $\mathcal{TL}$-category, and let $\De_{m+n}(r)$ be the cell module of $TL_{n}$ corresponding to index $r$.
Then we call~ $_{TL_{m}\o TL_{n}}\De_{m+n}(r)$, the restriction of $\De_{m+n}(r)$ as a ${TL_{m}\o TL_{n}}$-module, an $(m|n,r)$-walled  module as well as call
an $(m+n, r)$-cap diagram in $_{TL_{m}\o TL_{n}}\De_{m+n}(r)$ an  $(m|n, r)$-walled cap diagram. 
\end{defin}

For an $(m|n, r)$-walled cap diagram, we call the arcs crossing the wall \textit{through strings} and call the caps on the left  side of the wall \textit{left caps}, the caps on the right  side of the wall \textit{right caps}.
Furthermore, for an $(m|n, r)$-walled cap diagram, we can correspond to a triple $(s,l_{m},l_{n})$, where   $s$ is\textit{ the number of the through strings}, $l_{m}$ and
$l_{n}$ are \textit{the number of the left caps and the right caps}, respectively.  It is clear that $s+l_{m}+l_{n}=r$.

For instance, Figure \ref{fig2}(2) show that an $(11, 4)$-cap diagram in $\De_{11}(4)$ can be viewed  as
a $(6|5, 4)$-walled cap diagram in $_{TL_{6}\o TL_{5}}\De_{m+n}(4)$. Here, $\{ 3,8\}$, $\{ 6,7\}$
are through strings, $\{ 4,5\}$ is a left cap, $\{ 10,11\}$ is a right cap, and hence it correspond to the triple $(2,1,1)$.

\begin{figure}
	\[
	\xy
	(-8,7.5)*{\scriptstyle\bullet}="1s"; 
	(-1,7.5)*{\scriptstyle\bullet}="2s"; 
	(6,7.5)*{\scriptstyle\bullet}="3s";
	(13,7.5)*{\scriptstyle\bullet}="4s";
	(20,7.5)*{\scriptstyle\bullet} ="5s";
	(27,7.5)*{\scriptstyle\bullet}="6s";
	(34,7.5)*{\scriptstyle\bullet}="7s";
	(41,7.5)*{\scriptstyle\bullet}="8s";
	(48,7.5)*{\scriptstyle\bullet}="9s";
	(55,7.5)*{\scriptstyle\bullet}="10s";
	(62,7.5)*{\scriptstyle\bullet}="11s";
	(-8,17.5)*{\scriptstyle\bullet}="1d"; 
	(-1,17.5)*{\scriptstyle\bullet}="2d"; 
	(6,17.5)*{\scriptstyle\bullet}="3d";
	(13,17.5)*{\scriptstyle\bullet}="4d";
	(20,17.5)*{\scriptstyle\bullet} ="5d";
	(27,17.5)*{\scriptstyle\bullet}="6d";
	(34,17.5)*{\scriptstyle\bullet}="7d";
	(41,17.5)*{\scriptstyle\bullet}="8d";
	(48,17.5)*{\scriptstyle\bullet}="9d";
	(55,17.5)*{\scriptstyle\bullet}="10d";
	(62,17.5)*{\scriptstyle\bullet}="11d";
	(-8,19.7)*{\scriptstyle 1};
	(-1,19.7)*{\scriptstyle 2};
	(6,19.7)*{\scriptstyle 3};
	(13,19.7)*{\scriptstyle 4};
	(20,19.7)*{\scriptstyle 5};
	(27,19.7)*{\scriptstyle 6};
	(34,19.7)*{\scriptstyle 7};
	(41,19.7)*{\scriptstyle 8};
	(48,19.7)*{\scriptstyle 9};
	(55,19.7)*{\scriptstyle 10};
	(62,19.7)*{\scriptstyle 11};
	(6,5)*{}="3";
	(13,5)*{}="4";
	(20,5)*{}="5";
	(27,5)*{}="6";
	(34,5)*{}="7";
	(41,5)*{}="8";
	(48,5)*{}="9";
	(55,5)*{}="10"; 
	(62,5)*{}="11"; 
	"1d"; "2d" **\crv{(-6,15) & (-3.5,15)};
	"3d"; "1s" **\dir{-};
	"4d"; "5d" **\crv{(15,15) & (17.5,15)};
	"6d"; "6s" **\dir{-};
	"2s"; "3s" **\crv{(1,9.5) & (3.5,9.5)};
	"4s"; "5s" **\crv{(15,9.5) & (17.5,9.5)};
	"7d"; "8d" **\crv{(36,15) & (38.5,15)};
	"9d"; "7s" **\dir{-};
	"10d"; "10s" **\dir{-};
	"11d"; "11s" **\dir{-};
	"8s"; "9s" **\crv{(43,9.5) & (45.5,9.5)};
	(30.5,22)*{}="6.5";
	(30.5,3)*{}="7.5";
	"6.5"; "7.5" **\dir{.};
	\endxy
	\]
	\caption*{(1)}
\end{figure}

\begin{figure}
\[
\xy
(-8,5)*{\scriptstyle\bullet}; 
(-1,5)*{\scriptstyle\bullet}; 
(6,5)*{\scriptstyle\bullet} **\dir{};
(13,5)*{\scriptstyle\bullet} ; **\dir{};
(20,5)*{\scriptstyle\bullet} **\dir{};
(27,5)*{\scriptstyle\bullet}; 
(34,5)*{\scriptstyle\bullet}; 
(41,5)*{\scriptstyle\bullet};
(48,5)*{\scriptstyle\bullet} **\dir{};
(55,5)*{\scriptstyle\bullet}; 
(62,5)*{\scriptstyle\bullet}; **\dir{};
(-8,5)*{}="1";
(-1,5)*{}="2";
(6,5)*{}="3";
(13,5)*{}="4";
(20,5)*{}="5";
(27,5)*{}="6";
(34,5)*{}="7";
(41,5)*{}="8";
(48,5)*{}="9";
(55,5)*{}="10"; 
(62,5)*{}="11"; 
(-8,6.9)*{\scriptstyle 1};
(-1,6.9)*{\scriptstyle 2};
(6,6.9)*{\scriptstyle 3};
(13,6.9)*{\scriptstyle 4};
(20,6.9)*{\scriptstyle 5};
(27,6.9)*{\scriptstyle 6};
(34,6.9)*{\scriptstyle 7};
(41,6.9)*{\scriptstyle 8};
(48,6.9)*{\scriptstyle 9};
(55,6.9)*{\scriptstyle 10};
(62,6.9)*{\scriptstyle 11};
(-8,2)*{}="1.5";
"1"; "1.5" **\dir{-};
(-1,2)*{}="2.5";
"2"; "2.5" **\dir{-};
"3"; "8" **\crv{(17,-1.8) & (28,-1.8)};
"4"; "5" **\crv{(16,1.5) & (18,1.5)};
"6"; "7" **\crv{(29,2) & (32,2)};
(48,2)*{}="9.5";
"9"; "9.5" **\dir{-};
"10"; "11" **\crv{(57,2) & (60,2)};
(30.5,10)*{}="6.5";
(30.5,-3)*{}="7.5";
"6.5"; "7.5" **\dir{.};
\endxy
\]
\caption*{(2)}
	\caption{walled diagrams}
\label{fig2}
\end{figure}

Moreover, denote by $I$ the index set  consisting of the triples $(s,l_{m},l_{n})$ of all the $(m|n, r)$-walled cap diagrams,  and the triples $(s,l_{m},l_{n})$ is ordered lexicographically with `$s$' using natural numbers ordering, $l_{m}$~and~$l_{n}$~using   inverse ordering of natural numbers.

Let $W_{m|n}(s,l_{m},l_{n})$ be the $K$-subspace of  $_{TL_{m}\o TL_{n}}\De_{m+n}(r)$ spanned by all 
 walled  $(m|n, r)$-walled cap diagrams with indices  no more than $(s,l_{m},l_{n})$.

We then get a  chain in terms of the total ordering on $I$
$$0\subset \cdots \subset W_{m|n}(s^\pr,l_{m}^\pr,l_{n}^\pr)\subset W_{m|n}(s,l_{m},l_{n})\subset \cdots \subset _{TL_{m}\o TL_{n}}\De_{m+n}(r).\eqno(a)$$

Now, we claim: 

\begin{prop}\label{keyp} Keep the notation as above. Then:
	
(1)	Set $\De_{m|n}(s,l_{m},l_{n}):= W_{m|n}(s,l_{m},l_{n})/ W_{m|n}(s^\pr,l_{m}^\pr,l_{n}^\pr)$. Then we have
$$\De_{m|n}(s,l_{m},l_{n})\cong \De_{m}(l_{m})\o \De_{n}(l_{n})$$ as a TL$_{m}\o$TL$_{n}$-module;
	
(2) The chain (a)
is a composition series of $_{TL_{m}\o TL_{n}}\De_{m+n}(r)$.

\end{prop}

To prove this,
we first show that $ W(s,l_{m},l_{n})$ is a submodule of $_{TL_{m}\o TL_{n}}\De_{m+n}(r)$. Secondly, we prove Proposition \ref{keyp}(1), and then Proposition \ref{keyp}(2) follows immediately from \ref{keyp}(1). 

The key to prove Proposition \ref{keyp} is the following 
observations.

\begin{figure}[h]
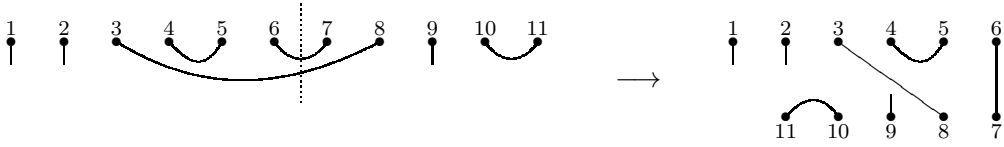

\[
\xy
(-8,5)*{\scriptstyle\bullet}; 
(-1,5)*{\scriptstyle\bullet}; 
(6,5)*{\scriptstyle\bullet} **\dir{};
(13,5)*{\scriptstyle\bullet} ; **\dir{};
(20,5)*{\scriptstyle\bullet} **\dir{};
(27,5)*{\scriptstyle\bullet}; 
(34,5)*{\scriptstyle\bullet}; 
(41,5)*{\scriptstyle\bullet};
(48,5)*{\scriptstyle\bullet} **\dir{};
(55,5)*{\scriptstyle\bullet}; 
(62,5)*{\scriptstyle\bullet}; **\dir{};
(-8,5)*{}="1";
(-1,5)*{}="2";
(6,5)*{}="3";
(13,5)*{}="4";
(20,5)*{}="5";
(27,5)*{}="6";
(34,5)*{}="7";
(41,5)*{}="8";
(48,5)*{}="9";
(55,5)*{}="10"; 
(62,5)*{}="11"; 
(-8,6.9)*{\scriptstyle 1};
(-1,6.9)*{\scriptstyle 2};
(6,6.9)*{\scriptstyle 3};
(13,6.9)*{\scriptstyle 4};
(20,6.9)*{\scriptstyle 5};
(27,6.9)*{\scriptstyle 6};
(34,6.9)*{\scriptstyle 7};
(41,6.9)*{\scriptstyle 8};
(48,6.9)*{\scriptstyle 9};
(55,6.9)*{\scriptstyle 10};
(62,6.9)*{\scriptstyle 11};
(-8,2)*{}="1.5";
"1"; "1.5" **\dir{-};
(-1,2)*{}="2.5";
"2"; "2.5" **\dir{-};
"3"; "8" **\crv{(17,-1.8) & (28,-1.8)};
"4"; "5" **\crv{(16,1.5) & (18,1.5)};
"6"; "7" **\crv{(29,2) & (32,2)};
(48,2)*{}="9.5";
"9"; "9.5" **\dir{-};
"10"; "11" **\crv{(57,2) & (60,2)};
(30.5,10)*{}="6.5";
(30.5,-3)*{}="7.5";
"6.5"; "7.5" **\dir{.};
\endxy
~~~~~~
\longrightarrow
~~~~~~
\xy
(-8,5)*{\scriptstyle\bullet}; 
(-1,5)*{\scriptstyle\bullet}; 
(6,5)*{\scriptstyle\bullet} **\dir{};
(13,5)*{\scriptstyle\bullet} ; **\dir{};
(20,5)*{\scriptstyle\bullet} **\dir{};
(27,5)*{\scriptstyle\bullet}; 
(27,-5)*{\scriptstyle\bullet}; 
(20,-5)*{\scriptstyle\bullet};
(13,-5)*{\scriptstyle\bullet} **\dir{};
(6,-5)*{\scriptstyle\bullet}; 
(-1,-5)*{\scriptstyle\bullet}; **\dir{};
(-8,5)*{}="1";
(-1,5)*{}="2";
(6,5)*{}="3";
(13,5)*{}="4";
(20,5)*{}="5";
(27,5)*{}="6";
(27,-5)*{}="7";
(20,-5)*{}="8";
(13,-5)*{}="9";
(6,-5)*{}="10"; 
(-1,-5)*{}="11"; 
(-8,6.9)*{\scriptstyle 1};
(-1,6.9)*{\scriptstyle 2};
(6,6.9)*{\scriptstyle 3};
(13,6.9)*{\scriptstyle 4};
(20,6.9)*{\scriptstyle 5};
(27,6.9)*{\scriptstyle 6};
(27,-6.9)*{\scriptstyle 7};
(20,-6.9)*{\scriptstyle 8};
(13,-6.9)*{\scriptstyle 9};
(6,-6.9)*{\scriptstyle 10};
(-1,-6.9)*{\scriptstyle 11};
(-8,2)*{}="1.5";
"1"; "1.5" **\dir{-};
(-1,2)*{}="2.5";
"2"; "2.5" **\dir{-};
"3"; "8" **\dir{-};
"4"; "5" **\crv{(16,1.5) & (18,1.5)};
"6"; "7" **\dir{-};
(13,-2)*{}="9.5";
"9"; "9.5" **\dir{-};
"11"; "10" **\crv{(1.5,-2) & (4,-2)};

\endxy
\]

	\caption{a $(6|5,4)$-walled cap diagram represented as a $(6,5)$-diagram}
\label{fig3}
\end{figure}

In fact, we notice that
an $(m|n, r)$-walled cap diagram with the triple $(s,l_{m},l_{n})$  also can be viewed as a so-called $(m,n)$-diagram, which is a diagram with $m$ vertices on the top  row and $n$ vertices on the bottom row labeling  the vertices in a
clockwise direction and the arc $\{i,j\}$ in the $(m,n)$-diagram is the same as in the $(m|n, r)$-walled cap diagram. (see Figure \ref{fig3} for an example.)

In Figure \ref{fig3}, the triple $(s,l_{m},l_{n})$ corresponding to the left diagram is $(2,1,1)$.
The triple in the right diagram can be read as follows: $s=2$ is the number of vertical arcs (therefore called the through strings), 
 $l_{m}=1$ is the number of the top horizontal arcs and $l_{n}=1$ is the number of the bottom horizontal arcs.

\begin{figure}[h]
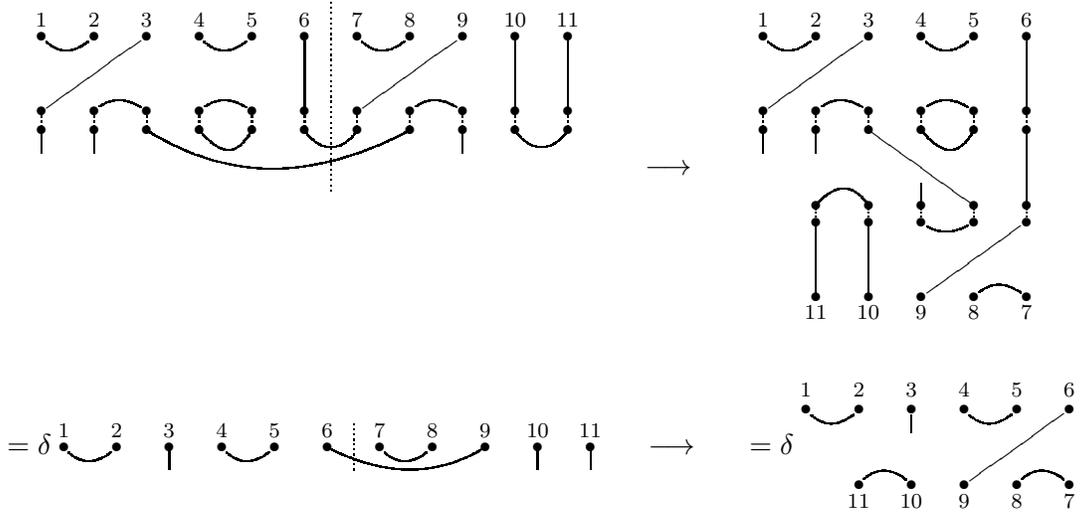

\[
\xy
(-8,7.5)*{\scriptstyle\bullet}="1s"; 
(-1,7.5)*{\scriptstyle\bullet}="2s"; 
(6,7.5)*{\scriptstyle\bullet}="3s";
(13,7.5)*{\scriptstyle\bullet}="4s";
(20,7.5)*{\scriptstyle\bullet} ="5s";
(27,7.5)*{\scriptstyle\bullet}="6s";
(34,7.5)*{\scriptstyle\bullet}="7s";
(41,7.5)*{\scriptstyle\bullet}="8s";
(48,7.5)*{\scriptstyle\bullet}="9s";
(55,7.5)*{\scriptstyle\bullet}="10s";
(62,7.5)*{\scriptstyle\bullet}="11s";
(-8,17.5)*{\scriptstyle\bullet}="1d"; 
(-1,17.5)*{\scriptstyle\bullet}="2d"; 
(6,17.5)*{\scriptstyle\bullet}="3d";
(13,17.5)*{\scriptstyle\bullet}="4d";
(20,17.5)*{\scriptstyle\bullet} ="5d";
(27,17.5)*{\scriptstyle\bullet}="6d";
(34,17.5)*{\scriptstyle\bullet}="7d";
(41,17.5)*{\scriptstyle\bullet}="8d";
(48,17.5)*{\scriptstyle\bullet}="9d";
(55,17.5)*{\scriptstyle\bullet}="10d";
(62,17.5)*{\scriptstyle\bullet}="11d";
(-8,19.7)*{\scriptstyle 1};
(-1,19.7)*{\scriptstyle 2};
(6,19.7)*{\scriptstyle 3};
(13,19.7)*{\scriptstyle 4};
(20,19.7)*{\scriptstyle 5};
(27,19.7)*{\scriptstyle 6};
(34,19.7)*{\scriptstyle 7};
(41,19.7)*{\scriptstyle 8};
(48,19.7)*{\scriptstyle 9};
(55,19.7)*{\scriptstyle 10};
(62,19.7)*{\scriptstyle 11};
(-8,5)*{\scriptstyle\bullet}; 
(-1,5)*{\scriptstyle\bullet}; 
(6,5)*{\scriptstyle\bullet} **\dir{};
(13,5)*{\scriptstyle\bullet} ; **\dir{};
(20,5)*{\scriptstyle\bullet} **\dir{};
(27,5)*{\scriptstyle\bullet}; 
(34,5)*{\scriptstyle\bullet}; 
(41,5)*{\scriptstyle\bullet};
(48,5)*{\scriptstyle\bullet} **\dir{};
(55,5)*{\scriptstyle\bullet}; 
(62,5)*{\scriptstyle\bullet}; **\dir{};
(-8,5)*{}="1";
(-1,5)*{}="2";
(6,5)*{}="3";
(13,5)*{}="4";
(20,5)*{}="5";
(27,5)*{}="6";
(34,5)*{}="7";
(41,5)*{}="8";
(48,5)*{}="9";
(55,5)*{}="10"; 
(62,5)*{}="11"; 
"1"; "1s" **\dir{.};
"2"; "2s" **\dir{.};
"3"; "3s" **\dir{.};
"4"; "4s" **\dir{.};
"5"; "5s" **\dir{.};
"6"; "6s" **\dir{.};
"7"; "7s" **\dir{.};
"8"; "8s" **\dir{.};
"9"; "9s" **\dir{.};
"10"; "10s" **\dir{.};
"11"; "11s" **\dir{.};
(-8,2)*{}="1.5";
"1"; "1.5" **\dir{-};
(-1,2)*{}="2.5";
"2"; "2.5" **\dir{-};
"3"; "8" **\crv{(17,-1.8) & (28,-1.8)};
"4"; "5" **\crv{(16,1.5) & (18,1.5)};
"6"; "7" **\crv{(29,2) & (32,2)};
(48,2)*{}="9.5";
"9"; "9.5" **\dir{-};
"10"; "11" **\crv{(57,2) & (60,2)};
(30.5,22)*{}="6.5";
(30.5,-3)*{}="7.5";
"6.5"; "7.5" **\dir{.};
"1d"; "2d" **\crv{(-6,15) & (-3.5,15)};
"3d"; "1s" **\dir{-};
"4d"; "5d" **\crv{(15,15) & (17.5,15)};
"6d"; "6s" **\dir{-};
"2s"; "3s" **\crv{(1,9.5) & (3.5,9.5)};
"4s"; "5s" **\crv{(15,9.5) & (17.5,9.5)};
"7d"; "8d" **\crv{(36,15) & (38.5,15)};
"9d"; "7s" **\dir{-};
"10d"; "10s" **\dir{-};
"11d"; "11s" **\dir{-};
"8s"; "9s" **\crv{(43,9.5) & (45.5,9.5)};
\endxy
~~~~~~
\longrightarrow
~~~~~~
\xy
(-8,5)*{\scriptstyle\bullet}; 
(-1,5)*{\scriptstyle\bullet}; 
(6,5)*{\scriptstyle\bullet} **\dir{};
(13,5)*{\scriptstyle\bullet} ; **\dir{};
(20,5)*{\scriptstyle\bullet} **\dir{};
(27,5)*{\scriptstyle\bullet}; 
(27,-5)*{\scriptstyle\bullet}; 
(20,-5)*{\scriptstyle\bullet};
(13,-5)*{\scriptstyle\bullet} **\dir{};
(6,-5)*{\scriptstyle\bullet}; 
(-1,-5)*{\scriptstyle\bullet}; **\dir{};
(-8,7.5)*{\scriptstyle\bullet}="1s"; 
(-1,7.5)*{\scriptstyle\bullet}="2s"; 
(6,7.5)*{\scriptstyle\bullet}="3s";
(13,7.5)*{\scriptstyle\bullet}="4s";
(20,7.5)*{\scriptstyle\bullet} ="5s";
(27,7.5)*{\scriptstyle\bullet}="6s";
(-8,17.5)*{\scriptstyle\bullet}="1d"; 
(-1,17.5)*{\scriptstyle\bullet}="2d"; 
(6,17.5)*{\scriptstyle\bullet}="3d";
(13,17.5)*{\scriptstyle\bullet}="4d";
(20,17.5)*{\scriptstyle\bullet} ="5d";
(27,17.5)*{\scriptstyle\bullet}="6d";
(-1,-7.2)*{\scriptstyle\bullet}="11s"; 
(6,-7.2)*{\scriptstyle\bullet}="10s"; 
(13,-7.2)*{\scriptstyle\bullet}="9s"; 
(20,-7.2)*{\scriptstyle\bullet}="8s"; 
(27,-7.2)*{\scriptstyle\bullet}="7s"; 
(-1,-17.2)*{\scriptstyle\bullet}="11d"; 
(6,-17.2)*{\scriptstyle\bullet}="10d"; 
(13,-17.2)*{\scriptstyle\bullet}="9d"; 
(20,-17.2)*{\scriptstyle\bullet}="8d"; 
(27,-17.2)*{\scriptstyle\bullet}="7d"; 
(-8,5)*{}="1";
(-1,5)*{}="2";
(6,5)*{}="3";
(13,5)*{}="4";
(20,5)*{}="5";
(27,5)*{}="6";
(27,-5)*{}="7";
(20,-5)*{}="8";
(13,-5)*{}="9";
(6,-5)*{}="10"; 
(-1,-5)*{}="11"; 
"7"; "7s" **\dir{.};
"8"; "8s" **\dir{.};
"9"; "9s" **\dir{.};
"10"; "10s" **\dir{.};
"11"; "11s" **\dir{.};
(-8,19.7)*{\scriptstyle 1};
(-1,19.7)*{\scriptstyle 2};
(6,19.7)*{\scriptstyle 3};
(13,19.7)*{\scriptstyle 4};
(20,19.7)*{\scriptstyle 5};
(27,19.7)*{\scriptstyle 6};
(27,-19.2)*{\scriptstyle 7};
(20,-19.2)*{\scriptstyle 8};
(13,-19.2)*{\scriptstyle 9};
(6,-19.2)*{\scriptstyle 10};
(-1,-19.2)*{\scriptstyle 11};
(-8,2)*{}="1.5";
"1"; "1.5" **\dir{-};
(-1,2)*{}="2.5";
"2"; "2.5" **\dir{-};
"3"; "8" **\dir{-};
"4"; "5" **\crv{(16,1.5) & (18,1.5)};
"6"; "7" **\dir{-};
(13,-2)*{}="9.5";
"9"; "9.5" **\dir{-};
"11"; "10" **\crv{(1.5,-2) & (4,-2)};
"1d"; "2d" **\crv{(-6,15) & (-3.5,15)};
"3d"; "1s" **\dir{-};
"4d"; "5d" **\crv{(15,15) & (17.5,15)};
"6d"; "6s" **\dir{-};
"2s"; "3s" **\crv{(1,9.5) & (3.5,9.5)};
"4s"; "5s" **\crv{(15,9.5) & (17.5,9.5)};
"1"; "1s" **\dir{.};
"2"; "2s" **\dir{.};
"3"; "3s" **\dir{.};
"4"; "4s" **\dir{.};
"5"; "5s" **\dir{.};
"6"; "6s" **\dir{.};
"11d"; "11s" **\dir{-};
"10d"; "10s" **\dir{-};
"9s"; "8s" **\crv{(15,-8.9) & (17.5,-8.9)};
"7s"; "9d" **\dir{-};
"8d"; "7d" **\crv{(22,-15) & (24.5,-15)};

\endxy
\]

\[
=\de~
\xy
(-8,0)*{\scriptstyle\bullet}="1d"; 
(-1,0)*{\scriptstyle\bullet}="2d"; 
(6,0)*{\scriptstyle\bullet}="3d";
(13,0)*{\scriptstyle\bullet}="4d";
(20,0)*{\scriptstyle\bullet} ="5d";
(27,0)*{\scriptstyle\bullet}="6d";
(34,0)*{\scriptstyle\bullet}="7d";
(41,0)*{\scriptstyle\bullet}="8d";
(48,0)*{\scriptstyle\bullet}="9d";
(55,0)*{\scriptstyle\bullet}="10d";
(62,0)*{\scriptstyle\bullet}="11d";
(-8,2.2)*{\scriptstyle 1};
(-1,2.2)*{\scriptstyle 2};
(6,2.2)*{\scriptstyle 3};
(13,2.2)*{\scriptstyle 4};
(20,2.2)*{\scriptstyle 5};
(27,2.2)*{\scriptstyle 6};
(34,2.2)*{\scriptstyle 7};
(41,2.2)*{\scriptstyle 8};
(48,2.2)*{\scriptstyle 9};
(55,2.2)*{\scriptstyle 10};
(62,2.2)*{\scriptstyle 11};
(-8,0)*{}="1";
(-1,0)*{}="2";
(6,0)*{}="3";
(13,0)*{}="4";
(20,0)*{}="5";
(27,0)*{}="6";
(34,0)*{}="7";
(41,0)*{}="8";
(48,0)*{}="9";
(55,0)*{}="10"; 
(62,0)*{}="11"; 
(30.5,3)*{}="6.5";
(30.5,-3)*{}="7.5";
"6.5"; "7.5" **\dir{.};
"1d"; "2d" **\crv{(-6,-2.5) & (-3.5,-2.5)};
"4d"; "5d" **\crv{(15,-2.5) & (17.5,-2.5)};
"7d"; "8d" **\crv{(36,-2.5) & (38.5,-2.5)};
(6,-3)*{}="3.5";
"3"; "3.5" **\dir{-};
(55,-3)*{}="10.5";
"10"; "10.5" **\dir{-};
(62,-3)*{}="11.5";
"11"; "11.5" **\dir{-};
"6"; "9" **\crv{(38.5,-6)};

\endxy
~~~~~~
\longrightarrow
~~~~~~
=\de~
\xy
(-8,5)*{\scriptstyle\bullet}="1";
(-1,5)*{\scriptstyle\bullet}="2";
(6,5)*{\scriptstyle\bullet} ="3";
(13,5)*{\scriptstyle\bullet} ="4";
(20,5)*{\scriptstyle\bullet}="5";
(27,5)*{\scriptstyle\bullet}="6";
(27,-5)*{\scriptstyle\bullet}="7"; 
(20,-5)*{\scriptstyle\bullet}="8";
(13,-5)*{\scriptstyle\bullet}="9";
(6,-5)*{\scriptstyle\bullet}="10";
(-1,-5)*{\scriptstyle\bullet}="11";
(-8,7.7)*{\scriptstyle 1};
(-1,7.7)*{\scriptstyle 2};
(6,7.7)*{\scriptstyle 3};
(13,7.7)*{\scriptstyle 4};
(20,7.7)*{\scriptstyle 5};
(27,7.7)*{\scriptstyle 6};
(27,-7.2)*{\scriptstyle 7};
(20,-7.2)*{\scriptstyle 8};
(13,-7.2)*{\scriptstyle 9};
(6,-7.2)*{\scriptstyle 10};
(-1,-7.2)*{\scriptstyle 11};
"1"; "2" **\crv{(-6,2.5) & (-3.5,2.5)};
(6,2)*{}="3.5";
"3"; "3.5" **\dir{-};
"4"; "5" **\crv{(15,2.5) & (17.5,2.5)};
(6,2)*{}="3.5";
"6"; "9" **\dir{-};
"8"; "7" **\crv{(22,-2.5) & (24.5,-2.5)};
"11"; "10" **\crv{(1,-2.5) & (3.5,-2.5)};
\endxy
\]
	\caption{ the action as a bimodule on  a walled diagram}
\label{fig4}
\end{figure}

Another observation is that,
$TL_{m}\o TL_{n}$-action on an $(m|n, r)$-walled cap diagram can be viewed as a left $TL_{m}$-right $TL_{n}^{op}$-action on  the corresponding $(m,n)$-diagram. Here, an $n$-diagram in $TL_{n}^{op}$ are obtained by rotating the $n$-diagram in $TL_{n}$ by a half turn. Figure \ref{fig4} is an example. 

With the above observations and by the properties of the concatenation of two Temperley-Lieb diagrams \cite{RSA}, the following lemma holds  immidiately. 

\begin{lemma}\label{updown}
Let $x$ be an  $(m|n, r)$-walled cap diagram with the triple $(s,l_{m},l_{n})$, and let $d$ be an $(m|n)$-diagrams in $TL_{m}\o TL_{n}$. Suppose the result diagram of the $TL_{m}\o TL_{n}$-action $d\cdot x$ has the triple $(s^\pr,l_{m}^\pr,l_{n}^\pr)$. Then we have $s^\pr \leqslant s, l_{m}\geqslant l_{m}^\pr$~and~$l_{n}\geqslant l_{n}^\pr$, that is, the number of the through
strings never  increases and the number of the caps never decreases.
\end{lemma}$\qed$

Due to above lemma, 
$ W_{m|n}(s,l_{m},l_{n})$ is a $TL_{m}\o TL_{n}$-module. Hence, the chain(a) is a  chain of modules.

Before proving Proposition \ref{keyp}(1), we also need the following result.

\begin{lemma}\label{key}
Let $u$ be an $(m,l_{m})$-cap diagram in $\De_{m}(l_{m}) $, and let $v$ be an $(n,l_{n})$-diagram in $\De_{n}(l_{n}) $ as well as
let $u\o v$, an $(m|n,l_{m}+l_{n})$-walled cap diagram, be the juxtaposition of $u$~and~$v$.  Suppose one can add   $s$  through strings on $u\o v$.
Then there is a unique $(m|n,l_{m}+l_{n}+s)$-walled cap diagram obtained by adding $s$  through strings on $u\o v$.
\end{lemma}
\proof    The result diagram by adding $s$  through strings on $u\o v$ can be obtained by the following steps (see Figure \ref{fig5} below). Firstly, all caps  are removed  directly from the diagram $u\o v$, and thus this leaves exactly vertices, say $i_{1}<i_{2}< \cdots<i_{m+n-2l_{m}}<i_{1}^{\pr}< i_{2}^\pr<\cdots<i_{m+n-2l_{n}}^\pr$, where $i_{j}$ is on the left of the wall and $i_{k}^{\pr}$ is on the right of the wall. It is an $(m-l_{m}|n-l_{n},0)$-walled cap diagram.

Secondly, we add $s$  through strings on the
the $(m-l_{m}|n-l_{n},0)$-walled cap diagram.
Because the through strings cross the wall and the $(m-l_{m}|n-l_{n},0)$-walled cap diagram is a half-diagram of an $(m+n-l_{m}-l_{n})$-diagram.  This means, for any added through string
$\{i_{j},i^{\pr}_{k}\}$, there is no single vertex $q$ with $i_{j}<q<i^{\pr}_{k}$ and moreover, for any two added through strings
$\{i_{j_{1}},i_{k_{1}}^{\pr}\}$, $\{i_{j_{2}},i_{k_{2}}^{\pr}\}$ with $i_{j_{1}}<i_{j_{2}}$, since $i_{j_{2}}<i_{k_{1}}^{\pr}$ and the strings cannot cross each other in a diagram, this implies  $i_{j_{1}}<i_{j_{2}}<i_{k_{2}}^{\pr}<i_{k_{1}}^{\pr}$.

Therefore, the $s$  through strings must be $\{i_{m+n-2l_{m}}, i_{1}^{\pr}\}, \{i_{m+n-2l_{m}-1},i_{2}^{\pr}\},$

$\cdots, \{i_{m+n-2l_{m}-s},i_{s}^{\pr}\}$. 

Finally, we get the result diagram by recovering all caps which we removed.
$\qed$

\begin{figure}[h]
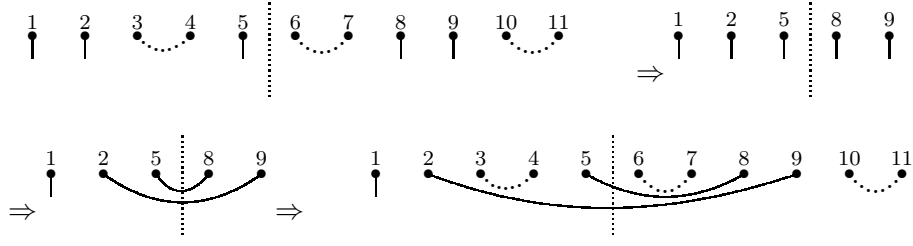

\[
\xy
(-8,5)*{\scriptstyle\bullet}; 
(-1,5)*{\scriptstyle\bullet}; 
(6,5)*{\scriptstyle\bullet} **\dir{};
(13,5)*{\scriptstyle\bullet} ; **\dir{};
(20,5)*{\scriptstyle\bullet} **\dir{};
(27,5)*{\scriptstyle\bullet}; 
(34,5)*{\scriptstyle\bullet}; 
(41,5)*{\scriptstyle\bullet};
(48,5)*{\scriptstyle\bullet} **\dir{};
(55,5)*{\scriptstyle\bullet}; 
(62,5)*{\scriptstyle\bullet}; **\dir{};
(-8,5)*{}="1";
(-1,5)*{}="2";
(6,5)*{}="3";
(13,5)*{}="4";
(20,5)*{}="5";
(27,5)*{}="6";
(34,5)*{}="7";
(41,5)*{}="8";
(48,5)*{}="9";
(55,5)*{}="10"; 
(62,5)*{}="11"; 
(-8,6.9)*{\scriptstyle 1};
(-1,6.9)*{\scriptstyle 2};
(6,6.9)*{\scriptstyle 3};
(13,6.9)*{\scriptstyle 4};
(20,6.9)*{\scriptstyle 5};
(27,6.9)*{\scriptstyle 6};
(34,6.9)*{\scriptstyle 7};
(41,6.9)*{\scriptstyle 8};
(48,6.9)*{\scriptstyle 9};
(55,6.9)*{\scriptstyle 10};
(62,6.9)*{\scriptstyle 11};
(-8,2)*{}="1.5";
"1"; "1.5" **\dir{-};
(-1,2)*{}="2.5";
"2"; "2.5" **\dir{-};
"3"; "4" **\crv{~*=<3pt>{.}(8,2.5) & (10.5,2.5)};
(20,2)*{}="5.5";
"5"; "5.5" **\dir{-};
"6"; "7" **\crv{~*=<3pt>{.}(29,2) & (32,2)};
(41,2)*{}="8.5";
"8"; "8.5" **\dir{-};
(48,2)*{}="9.5";
"9"; "9.5" **\dir{-};
"10"; "11" **\crv{~*=<3pt>{.}(57,2) & (60,2)};
(23.5,10)*{}="5.5";
(23.5,-3)*{}="6.5";
"5.5"; "6.5" **\dir{.};
\endxy
~~~~~~
\Rightarrow
\xy
(-8,5)*{\scriptstyle\bullet}; 
(-1,5)*{\scriptstyle\bullet}; 
(6,5)*{\scriptstyle\bullet} **\dir{};
(13,5)*{\scriptstyle\bullet} **\dir{};
(20,5)*{\scriptstyle\bullet} **\dir{};
(-8,5)*{}="1";
(-1,5)*{}="2";
(6,5)*{}="5";
(13,5)*{}="8";
(20,5)*{}="9";
(-8,7.2)*{\scriptstyle 1};
(-1,7.2)*{\scriptstyle 2};
(6,7.2)*{\scriptstyle 5};
(13,7.2)*{\scriptstyle 8};
(20,7.2)*{\scriptstyle 9};
(-8,2)*{}="1.5";
"1"; "1.5" **\dir{-};
(-1,2)*{}="2.5";
"2"; "2.5" **\dir{-};
(6,2)*{}="5.5";
"5"; "5.5" **\dir{-};
(13,2)*{}="8.5";
"8"; "8.5" **\dir{-};
(20,2)*{}="9.5";
"9"; "9.5" **\dir{-};
(9.5,10)*{}="5.5";
(9.5,-3)*{}="8.5";
"5.5"; "8.5" **\dir{.};
\endxy
\]

\[
\Rightarrow
\xy
(-8,5)*{\scriptstyle\bullet}; 
(-1,5)*{\scriptstyle\bullet}; 
(6,5)*{\scriptstyle\bullet} **\dir{};
(13,5)*{\scriptstyle\bullet} **\dir{};
(20,5)*{\scriptstyle\bullet} **\dir{};
(-8,5)*{}="1";
(-1,5)*{}="2";
(6,5)*{}="5";
(13,5)*{}="8";
(20,5)*{}="9";
(-8,7.2)*{\scriptstyle 1};
(-1,7.2)*{\scriptstyle 2};
(6,7.2)*{\scriptstyle 5};
(13,7.2)*{\scriptstyle 8};
(20,7.2)*{\scriptstyle 9};
(-8,2)*{}="1.5";
"1"; "1.5" **\dir{-};
(9.5,10)*{}="5.5";
(9.5,-3)*{}="8.5";
"5.5"; "8.5" **\dir{.};
"5"; "8" **\crv{(8,2) & (10.5,2)};
"2"; "9" **\crv{(6,0) & (13,0)};

\endxy
\Rightarrow
~~~~~~
\xy
(-8,5)*{\scriptstyle\bullet}; 
(-1,5)*{\scriptstyle\bullet}; 
(6,5)*{\scriptstyle\bullet} **\dir{};
(13,5)*{\scriptstyle\bullet} ; **\dir{};
(20,5)*{\scriptstyle\bullet} **\dir{};
(27,5)*{\scriptstyle\bullet}; 
(34,5)*{\scriptstyle\bullet}; 
(41,5)*{\scriptstyle\bullet};
(48,5)*{\scriptstyle\bullet} **\dir{};
(55,5)*{\scriptstyle\bullet}; 
(62,5)*{\scriptstyle\bullet}; **\dir{};
(-8,5)*{}="1";
(-1,5)*{}="2";
(6,5)*{}="3";
(13,5)*{}="4";
(20,5)*{}="5";
(27,5)*{}="6";
(34,5)*{}="7";
(41,5)*{}="8";
(48,5)*{}="9";
(55,5)*{}="10"; 
(62,5)*{}="11"; 
(-8,7.2)*{\scriptstyle 1};
(-1,7.2)*{\scriptstyle 2};
(6,7.2)*{\scriptstyle 3};
(13,7.2)*{\scriptstyle 4};
(20,7.2)*{\scriptstyle 5};
(27,7.2)*{\scriptstyle 6};
(34,7.2)*{\scriptstyle 7};
(41,7.2)*{\scriptstyle 8};
(48,7.2)*{\scriptstyle 9};
(55,7.2)*{\scriptstyle 10};
(62,7.2)*{\scriptstyle 11};
(-8,2)*{}="1.5";
"1"; "1.5" **\dir{-};
"3"; "4" **\crv{~*=<3pt>{.}(8,2.5) & (10.5,2.5)};
"5"; "8" **\crv{(27,1) & (34,1)};
"6"; "7" **\crv{~*=<3pt>{.}(29,2) & (32,2)};
"2"; "9" **\crv{(15,-1) & (30,-1)};
"10"; "11" **\crv{~*=<3pt>{.}(57,2) & (60,2)};
(23.5,10)*{}="5.5";
(23.5,-3)*{}="6.5";
"5.5"; "6.5" **\dir{.};
\endxy
\]

	\caption{ adding  through strings on  a walled cap diagram}
\label{fig5}
\end{figure}

\

We are now in the position to give a proof of Proposition \ref{keyp}(1).

\proof
Since a $TL_{m}\o TL_{n}$-module can be viewed as   a left $TL_{m}$-right $TL_{n}^{op}$-bimodule, it is sufficient to prove $\De_{m|n}(s,l_{m},l_{n})\cong \De_{m}(l_{m})\o \De_{n}(l_{n})$ as a $TL_{m}$-module.

Firstly, we give an one-to-one map from $\De_{m}(l_{m})\o \De_{n}(l_{n})$ to  $\De_{m|n}(s,l_{m},l_{n})$. Due to Lemma \ref{key}, given an $(m|n,l_{m}+l_{n})$-cap diagram $u\o v$ in $\De_{m}(l_{m})\o \De_{n}(l_{n})$, there is a unique diagram obtained by adding $s$  through strings on $u\o v$ in $\De_{m|n}(s,l_{m},l_{n})$. Consequently, let such diagram correspond to $u\o v$. we then denote this map by $\sigma$ and extend linearly.

We next show that $\sigma$ preserves the left $TL_{m}$-action.

Precisely, let $u=C_{S_{1}}^{l_{m}}$ and $v=C_{S_{2}}^{l_{n}}$ be two cap diagrams in $\De_{m}(l_{m})$ and $ \De_{n}(l_{n})$ respectively. For any $m$-diagram $C_{S,T}^{l}$ in $TL_{m}$(since all $m$-diagram form a cellular basis),
if $l> l_{m}$, then $C_{S,T}^{l}\cdot (C_{S_{1}}^{l_{m}} \o C_{S_{2}}^{l_{n}})= (C_{S,T}^{l}\cdot C_{S_{1}}^{l_{m}} )\o C_{S_{2}}^{l_{n}}=0$. On the other hand, the number of the left cups of walled cap diagram $C_{S,T}^{l}\cdot \sigma (C_{S_{1}}^{l_{m}} \o C_{S_{2}}^{l_{n}})$ is increases and the number of the through strings is never increase due to Lemma \ref{updown}, hence it is equal to $0$ in $\De_{m|n}(s,l_{m},l_{n})$.

If $l\leqslant l_{m}$, due to previous observation in Figure \ref{fig4},  the diagram $C_{S,T}^{l}\cdot \sigma (C_{S_{1}}^{l_{m}} \o C_{S_{2}}^{l_{n}})$ can be obtained by  following steps(see Figure \ref{fig6} below): Firstly, cutting the diagram  $C_{S,T}^{l}\cdot \sigma (C_{S_{1}}^{l_{m}} \o C_{S_{2}}^{l_{n}})$ along the wall, the diagram becomes   the juxtaposition of $C_{S,T}^{l}\cdot C_{S_{1}}^{l_{m}}  $ and $C_{S_{2}}^{l_{n}}$, say $ (C_{S,T}^{l}\cdot C_{S_{1}}^{l_{m}} )\o C_{S_{2}}^{l_{n}}$.  It remains to  add $s$ through strings on $ (C_{S,T}^{l}\cdot C_{S_{1}}^{l_{m}} )\o C_{S_{2}}^{l_{n}}$. It is, however, just $\sigma(C_{S,T}^{l}\cdot (C_{S_{1}}^{l_{m}} \o C_{S_{2}}^{l_{n}}))$ by Lemma \ref{key}.


$\qed$

\begin{figure}[h]
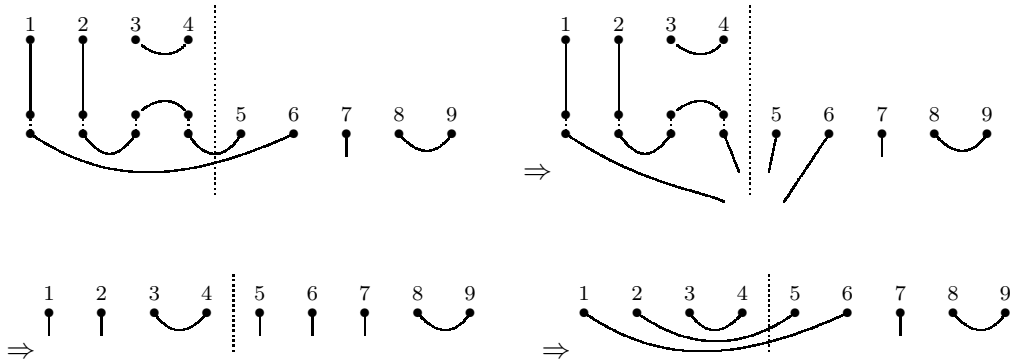


\[
\xy
(-8,7.5)*{\scriptstyle\bullet}="1s"; 
(-1,7.5)*{\scriptstyle\bullet}="2s"; 
(6,7.5)*{\scriptstyle\bullet}="3s";
(13,7.5)*{\scriptstyle\bullet}="4s";
(-8,17.5)*{\scriptstyle\bullet}="1d"; 
(-1,17.5)*{\scriptstyle\bullet}="2d"; 
(6,17.5)*{\scriptstyle\bullet}="3d";
(13,17.5)*{\scriptstyle\bullet}="4d";
(-8,19.7)*{\scriptstyle 1};
(-1,19.7)*{\scriptstyle 2};
(6,19.7)*{\scriptstyle 3};
(13,19.7)*{\scriptstyle 4};
(20,7.7)*{\scriptstyle 5};
(27,7.7)*{\scriptstyle 6};
(34,7.7)*{\scriptstyle 7};
(41,7.7)*{\scriptstyle 8};
(48,7.7)*{\scriptstyle 9};
(-8,5)*{\scriptstyle\bullet}; 
(-1,5)*{\scriptstyle\bullet}; 
(6,5)*{\scriptstyle\bullet} **\dir{};
(13,5)*{\scriptstyle\bullet} ; **\dir{};
(20,5)*{\scriptstyle\bullet} **\dir{};
(27,5)*{\scriptstyle\bullet}; 
(34,5)*{\scriptstyle\bullet}; 
(41,5)*{\scriptstyle\bullet};
(48,5)*{\scriptstyle\bullet} **\dir{};
(-8,5)*{}="1";
(-1,5)*{}="2";
(6,5)*{}="3";
(13,5)*{}="4";
(20,5)*{}="5";
(27,5)*{}="6";
(34,5)*{}="7";
(41,5)*{}="8";
(48,5)*{}="9";
"1"; "1s" **\dir{.};
"2"; "2s" **\dir{.};
"3"; "3s" **\dir{.};
"4"; "4s" **\dir{.};
(34,2)*{}="7.5";
"7"; "7.5" **\dir{-};
"1"; "6" **\crv{(2,-1.8) & (12,-1.8)};
"2"; "3" **\crv{(1.5,1.5) & (3.5,1.5)};
"4"; "5" **\crv{(15,1.5) & (18,1.5)};
"8"; "9" **\crv{(43.5,2) & (46,2)};
(16.5,22)*{}="4.5";
(16.5,-3)*{}="5.5";
"4.5"; "5.5" **\dir{.};
"1d"; "1s" **\dir{-};
"2d"; "2s" **\dir{-};
"3d"; "4d" **\crv{(8.5,15) & (11.5,15)};
"3s"; "4s" **\crv{(8.5,10) & (11.5,10)};
\endxy
~~~~~~
\Rightarrow
\xy
(-8,7.5)*{\scriptstyle\bullet}="1s"; 
(-1,7.5)*{\scriptstyle\bullet}="2s"; 
(6,7.5)*{\scriptstyle\bullet}="3s";
(13,7.5)*{\scriptstyle\bullet}="4s";
(-8,17.5)*{\scriptstyle\bullet}="1d"; 
(-1,17.5)*{\scriptstyle\bullet}="2d"; 
(6,17.5)*{\scriptstyle\bullet}="3d";
(13,17.5)*{\scriptstyle\bullet}="4d";
(-8,19.7)*{\scriptstyle 1};
(-1,19.7)*{\scriptstyle 2};
(6,19.7)*{\scriptstyle 3};
(13,19.7)*{\scriptstyle 4};
(20,7.7)*{\scriptstyle 5};
(27,7.7)*{\scriptstyle 6};
(34,7.7)*{\scriptstyle 7};
(41,7.7)*{\scriptstyle 8};
(48,7.7)*{\scriptstyle 9};
(-8,5)*{\scriptstyle\bullet}; 
(-1,5)*{\scriptstyle\bullet}; 
(6,5)*{\scriptstyle\bullet} **\dir{};
(13,5)*{\scriptstyle\bullet} ; **\dir{};
(20,5)*{\scriptstyle\bullet} **\dir{};
(27,5)*{\scriptstyle\bullet}; 
(34,5)*{\scriptstyle\bullet}; 
(41,5)*{\scriptstyle\bullet};
(48,5)*{\scriptstyle\bullet} **\dir{};
(-8,5)*{}="1";
(-1,5)*{}="2";
(6,5)*{}="3";
(13,5)*{}="4";
(20,5)*{}="5";
(27,5)*{}="6";
(34,5)*{}="7";
(41,5)*{}="8";
(48,5)*{}="9";
"1"; "1s" **\dir{.};
"2"; "2s" **\dir{.};
"3"; "3s" **\dir{.};
"4"; "4s" **\dir{.};
(34,2)*{}="7.5";
"7"; "7.5" **\dir{-};
(13,-4)*{}="3.5";
"1"; "3.5" **\crv{(2,-1.8) & (10,-2.5)};
(21,-4)*{}="5.5";
"5.5"; "6" **\crv{(21,-4) & (21,-4)};
"2"; "3" **\crv{(1.5,1.5) & (3.5,1.5)};
(15,0)*{}="4.5";
"4"; "4.5" **\crv{(15,0) & (15,0)};
(19,0)*{}="4.6";
"5"; "4.6" **\crv{(19,0) & (19,0)};
"8"; "9" **\crv{(43.5,2) & (46,2)};
(16.5,22)*{}="4.5";
(16.5,-3)*{}="5.5";
"4.5"; "5.5" **\dir{.};
"1d"; "1s" **\dir{-};
"2d"; "2s" **\dir{-};
"3d"; "4d" **\crv{(8.5,15) & (11.5,15)};
"3s"; "4s" **\crv{(8.5,10) & (11.5,10)};
\endxy
\]

\[
\Rightarrow
\xy
(-8,7.7)*{\scriptstyle 1};
(-1,7.7)*{\scriptstyle 2};
(6,7.7)*{\scriptstyle 3};
(13,7.7)*{\scriptstyle 4};
(20,7.7)*{\scriptstyle 5};
(27,7.7)*{\scriptstyle 6};
(34,7.7)*{\scriptstyle 7};
(41,7.7)*{\scriptstyle 8};
(48,7.7)*{\scriptstyle 9};
(-8,5)*{\scriptstyle\bullet}; 
(-1,5)*{\scriptstyle\bullet}; 
(6,5)*{\scriptstyle\bullet} **\dir{};
(13,5)*{\scriptstyle\bullet} ; **\dir{};
(20,5)*{\scriptstyle\bullet} **\dir{};
(27,5)*{\scriptstyle\bullet}; 
(34,5)*{\scriptstyle\bullet}; 
(41,5)*{\scriptstyle\bullet};
(48,5)*{\scriptstyle\bullet} **\dir{};
(-8,5)*{}="1";
(-1,5)*{}="2";
(6,5)*{}="3";
(13,5)*{}="4";
(20,5)*{}="5";
(27,5)*{}="6";
(34,5)*{}="7";
(41,5)*{}="8";
(48,5)*{}="9";
(-8,2)*{}="1.5";
"1"; "1.5" **\dir{-};
(-1,2)*{}="2.5";
"2"; "2.5" **\dir{-};
(20,2)*{}="5.5";
"5"; "5.5" **\dir{-};
(27,2)*{}="6.5";
"6"; "6.5" **\dir{-};
(34,2)*{}="7.5";
"7"; "7.5" **\dir{-};
"3"; "4" **\crv{(8,2) & (10.5,2)};
"8"; "9" **\crv{(43.5,2) & (46,2)};
(16.5,10)*{}="4.5";
(16.5,0)*{}="5.5";
"4.5"; "5.5" **\dir{.};
\endxy
~~~~~~
\Rightarrow
\xy
(-8,7.7)*{\scriptstyle 1};
(-1,7.7)*{\scriptstyle 2};
(6,7.7)*{\scriptstyle 3};
(13,7.7)*{\scriptstyle 4};
(20,7.7)*{\scriptstyle 5};
(27,7.7)*{\scriptstyle 6};
(34,7.7)*{\scriptstyle 7};
(41,7.7)*{\scriptstyle 8};
(48,7.7)*{\scriptstyle 9};
(-8,5)*{\scriptstyle\bullet}; 
(-1,5)*{\scriptstyle\bullet}; 
(6,5)*{\scriptstyle\bullet} **\dir{};
(13,5)*{\scriptstyle\bullet} ; **\dir{};
(20,5)*{\scriptstyle\bullet} **\dir{};
(27,5)*{\scriptstyle\bullet}; 
(34,5)*{\scriptstyle\bullet}; 
(41,5)*{\scriptstyle\bullet};
(48,5)*{\scriptstyle\bullet} **\dir{};
(-8,5)*{}="1";
(-1,5)*{}="2";
(6,5)*{}="3";
(13,5)*{}="4";
(20,5)*{}="5";
(27,5)*{}="6";
(34,5)*{}="7";
(41,5)*{}="8";
(48,5)*{}="9";
(34,2)*{}="7.5";
"7"; "7.5" **\dir{-};
"3"; "4" **\crv{(8,2) & (10.5,2)};
"8"; "9" **\crv{(43.5,2) & (46,2)};
(16.5,10)*{}="4.5";
(16.5,0)*{}="5.5";
"4.5"; "5.5" **\dir{.};
"1"; "6" **\crv{(2,-1.8) & (12,-1.8)};
"2"; "5" **\crv{(6,0) & (13,0)};
\endxy
\]

	\caption{ the left $TL_{4}$-action}
\label{fig6}
\end{figure}

Due to Proposition \ref{keyp}(2) and equality (6), the structure constant $a_{(p, q, r)}^{(m|n)}$ in Theorem \ref{cons} is the multiplicity of the simple module $\De_{m}(p)\o \De_{n}(q)$ occurred as a composition factor in the chain (a), a composition series of $_{TL_{m}\o TL_{n}}\De_{m+n}(r)$.

Moreover, by Proposition \ref{keyp}(1), suppose
$\De_{m|n}(s,p,q)$ ($\cong \De_{m}(p)\o \De_{n}(q)$) is a composition factor in the chain (a). Then the triple $(s,p,q)$ must be satisfied:

\[
\left\{
\begin{array}{ll}
s+p+q=r  ~~~(total~r ~arcs~),\\
m-s\geqslant 2p  ~~~(p ~left~ caps),\\
n-s\geqslant 2q  ~~~(q ~right~ caps).
\end{array}
\right.
\]

Conversely, if $\De_{m}(p)\o \De_{n}(q)$ occur as a composition factor in the chain (a), then $\De_{m|n}(s,p,q)$, with $s=r-p-q$, is the unique $(m|n,r)$-walled cap module isomorphic to $ \De_{m}(p)\o \De_{n}(q)$. Therefore, the multiplicity equals $1$.

Consequently, due to above facts, Theorem \ref{cons} follows.

Now, we end this paper by a simple examples to illustrate our results.

\textbf{Example 2} ~~For $_{TL_{4}\o TL_{3}}\De_{7}(2)$, we display a composition series as follows:

$0\subset W_{4|3}(0,2,0)\subset W_{4|3}(0,1,1)\subset\\~~~~~~~~~~~~~~~~~~~~~~~~~~~~~~~~~~~ W_{4|3}(1,1,0)\subset W_{4|3}(1,0,1)\subset W_{4|3}(2,0,0) = _{TL_{4}\o TL_{3}}\De_{7}(2).$

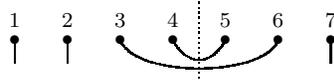
\begin{figure}[h]
	\[
	\xy
	(-8,7.7)*{\scriptstyle 1};
	(-1,7.7)*{\scriptstyle 2};
	(6,7.7)*{\scriptstyle 3};
	(13,7.7)*{\scriptstyle 4};
	(20,7.7)*{\scriptstyle 5};
	(27,7.7)*{\scriptstyle 6};
	(34,7.7)*{\scriptstyle 7};
	(-8,5)*{\scriptstyle\bullet}; 
	(-1,5)*{\scriptstyle\bullet}; 
	(6,5)*{\scriptstyle\bullet} **\dir{};
	(13,5)*{\scriptstyle\bullet} ; **\dir{};
	(20,5)*{\scriptstyle\bullet} **\dir{};
	(27,5)*{\scriptstyle\bullet}; 
	(34,5)*{\scriptstyle\bullet}; 
	(-8,5)*{}="1";
	(-1,5)*{}="2";
	(6,5)*{}="3";
	(13,5)*{}="4";
	(20,5)*{}="5";
	(27,5)*{}="6";
	(34,5)*{}="7";
	(-8,2)*{}="1.5";
	"1"; "1.5" **\dir{-};
	(-1,2)*{}="2.5";
	"2"; "2.5" **\dir{-};
	(6,2)*{}="3.5";
		(27,2)*{}="6.5";
		(34,2)*{}="7.5";
		"7"; "7.5" **\dir{-};
	"3"; "6" **\crv{(8,0) & (25,0)};
	"4"; "5" **\crv{(15,1.5) & (18,1.5)};
	(16.5,10)*{}="4.5";
	(16.5,0)*{}="5.5";
	"4.5"; "5.5" **\dir{.};
	\endxy 
	\]	
	\caption*{ $\De_{4|3}(2,0,0)\cong \De_{4}(0)\o \De_{3}(0)$}
\end{figure}
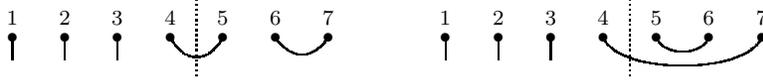
\begin{figure}[h]
	\[
	\xy
	(-8,7.7)*{\scriptstyle 1};
	(-1,7.7)*{\scriptstyle 2};
	(6,7.7)*{\scriptstyle 3};
	(13,7.7)*{\scriptstyle 4};
	(20,7.7)*{\scriptstyle 5};
	(27,7.7)*{\scriptstyle 6};
	(34,7.7)*{\scriptstyle 7};
	(-8,5)*{\scriptstyle\bullet}; 
	(-1,5)*{\scriptstyle\bullet}; 
	(6,5)*{\scriptstyle\bullet} **\dir{};
	(13,5)*{\scriptstyle\bullet} ; **\dir{};
	(20,5)*{\scriptstyle\bullet} **\dir{};
	(27,5)*{\scriptstyle\bullet}; 
	(34,5)*{\scriptstyle\bullet}; 
	(-8,5)*{}="1";
	(-1,5)*{}="2";
	(6,5)*{}="3";
	(13,5)*{}="4";
	(20,5)*{}="5";
	(27,5)*{}="6";
	(34,5)*{}="7";
	(-8,2)*{}="1.5";
	"1"; "1.5" **\dir{-};
	(-1,2)*{}="2.5";
	"2"; "2.5" **\dir{-};
	(6,2)*{}="3.5";
	"3"; "3.5" **\dir{-};
	"4"; "5" **\crv{(15,1.5) & (18,1.5)};
	"6"; "7" **\crv{(29,2) & (32,2)};
	(16.5,10)*{}="4.5";
	(16.5,0)*{}="5.5";
	"4.5"; "5.5" **\dir{.};
	\endxy 
~~~~~~~~~~~	\xy
(-8,7.7)*{\scriptstyle 1};
(-1,7.7)*{\scriptstyle 2};
(6,7.7)*{\scriptstyle 3};
(13,7.7)*{\scriptstyle 4};
(20,7.7)*{\scriptstyle 5};
(27,7.7)*{\scriptstyle 6};
(34,7.7)*{\scriptstyle 7};
(-8,5)*{\scriptstyle\bullet}; 
(-1,5)*{\scriptstyle\bullet}; 
(6,5)*{\scriptstyle\bullet} **\dir{};
(13,5)*{\scriptstyle\bullet} ; **\dir{};
(20,5)*{\scriptstyle\bullet} **\dir{};
(27,5)*{\scriptstyle\bullet}; 
(34,5)*{\scriptstyle\bullet}; 
(-8,5)*{}="1";
(-1,5)*{}="2";
(6,5)*{}="3";
(13,5)*{}="4";
(20,5)*{}="5";
(27,5)*{}="6";
(34,5)*{}="7";
(-8,2)*{}="1.5";
"1"; "1.5" **\dir{-};
(-1,2)*{}="2.5";
"2"; "2.5" **\dir{-};
(6,2)*{}="3.5";
"3"; "3.5" **\dir{-};
"4"; "7" **\crv{(15,0) & (32,0)};
"5"; "6" **\crv{(21,2.5) & (26,2.5)};
(16.5,10)*{}="4.5";
(16.5,0)*{}="5.5";
"4.5"; "5.5" **\dir{.};
\endxy 
	\]
	
	\caption*{ $\De_{4|3}(1,0,1)\cong \De_{4}(0)\o \De_{3}(1)$}
\end{figure}
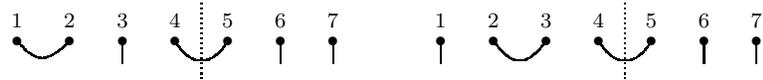
\begin{figure}[h]
		
	\[
	\xy
	(-8,7.7)*{\scriptstyle 1};
	(-1,7.7)*{\scriptstyle 2};
	(6,7.7)*{\scriptstyle 3};
	(13,7.7)*{\scriptstyle 4};
	(20,7.7)*{\scriptstyle 5};
	(27,7.7)*{\scriptstyle 6};
	(34,7.7)*{\scriptstyle 7};
	(-8,5)*{\scriptstyle\bullet}; 
	(-1,5)*{\scriptstyle\bullet}; 
	(6,5)*{\scriptstyle\bullet} **\dir{};
	(13,5)*{\scriptstyle\bullet} ; **\dir{};
	(20,5)*{\scriptstyle\bullet} **\dir{};
	(27,5)*{\scriptstyle\bullet}; 
	(34,5)*{\scriptstyle\bullet}; 
	(-8,5)*{}="1";
	(-1,5)*{}="2";
	(6,5)*{}="3";
	(13,5)*{}="4";
	(20,5)*{}="5";
	(27,5)*{}="6";
	(34,5)*{}="7";
(6,2)*{}="3.5";
"3"; "3.5" **\dir{-};
	(27,2)*{}="6.5";
	"6"; "6.5" **\dir{-};
	(34,2)*{}="7.5";
	"7"; "7.5" **\dir{-};
	"1"; "2" **\crv{(-5.5,2) & (-4,2)};
	"4"; "5" **\crv{(15,1.5) & (18,1.5)};
	(16.5,10)*{}="4.5";
	(16.5,0)*{}="5.5";
	"4.5"; "5.5" **\dir{.};
	\endxy ~~~~~~~~~~
	\xy
(-8,7.7)*{\scriptstyle 1};
(-1,7.7)*{\scriptstyle 2};
(6,7.7)*{\scriptstyle 3};
(13,7.7)*{\scriptstyle 4};
(20,7.7)*{\scriptstyle 5};
(27,7.7)*{\scriptstyle 6};
(34,7.7)*{\scriptstyle 7};
(-8,5)*{\scriptstyle\bullet}; 
(-1,5)*{\scriptstyle\bullet}; 
(6,5)*{\scriptstyle\bullet} **\dir{};
(13,5)*{\scriptstyle\bullet} ; **\dir{};
(20,5)*{\scriptstyle\bullet} **\dir{};
(27,5)*{\scriptstyle\bullet}; 
(34,5)*{\scriptstyle\bullet}; 
(-8,5)*{}="1";
(-1,5)*{}="2";
(6,5)*{}="3";
(13,5)*{}="4";
(20,5)*{}="5";
(27,5)*{}="6";
(34,5)*{}="7";
	(-8,2)*{}="1.5";
(6,2)*{}="3.5";
	"1"; "1.5" **\dir{-};
(27,2)*{}="6.5";
"6"; "6.5" **\dir{-};
(34,2)*{}="7.5";
"7"; "7.5" **\dir{-};
"2"; "3" **\crv{(1.5,1.5) & (3.5,1.5)};
"4"; "5" **\crv{(15,1.5) & (18,1.5)};
(16.5,10)*{}="4.5";
(16.5,0)*{}="5.5";
"4.5"; "5.5" **\dir{.};
\endxy 
\]
	\[
\xy
(-8,7.7)*{\scriptstyle 1};
(-1,7.7)*{\scriptstyle 2};
(6,7.7)*{\scriptstyle 3};
(13,7.7)*{\scriptstyle 4};
(20,7.7)*{\scriptstyle 5};
(27,7.7)*{\scriptstyle 6};
(34,7.7)*{\scriptstyle 7};
(-8,5)*{\scriptstyle\bullet}; 
(-1,5)*{\scriptstyle\bullet}; 
(6,5)*{\scriptstyle\bullet} **\dir{};
(13,5)*{\scriptstyle\bullet} ; **\dir{};
(20,5)*{\scriptstyle\bullet} **\dir{};
(27,5)*{\scriptstyle\bullet}; 
(34,5)*{\scriptstyle\bullet}; 
(-8,5)*{}="1";
(-1,5)*{}="2";
(6,5)*{}="3";
(13,5)*{}="4";
(20,5)*{}="5";
(27,5)*{}="6";
(34,5)*{}="7";
(-8,2)*{}="1.5";
"1"; "1.5" **\dir{-};
(-1,2)*{}="2.5";
(20,2)*{}="5.5";
(27,2)*{}="6.5";
"6"; "6.5" **\dir{-};
(34,2)*{}="7.5";
"7"; "7.5" **\dir{-};
"3"; "4" **\crv{(8,2) & (10.5,2)};
"2"; "5" **\crv{(1.5,0) & (17,0)};
(16.5,10)*{}="4.5";
(16.5,0)*{}="5.5";
"4.5"; "5.5" **\dir{.};
\endxy 
\]
	
	\caption*{ $\De_{4|3}(1,1,0)\cong \De_{4}(1)\o \De_{3}(0)$}
\end{figure}

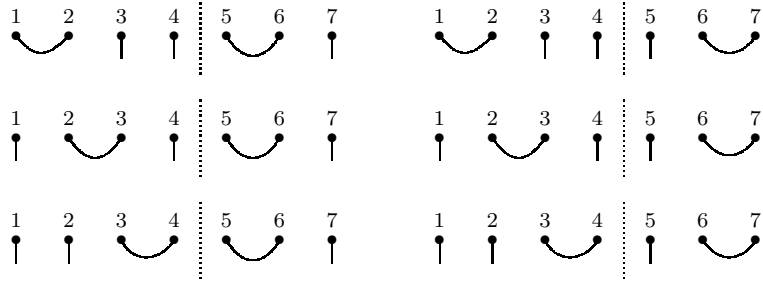
\begin{figure}[h!]
	\[
	\xy
	(-8,7.7)*{\scriptstyle 1};
	(-1,7.7)*{\scriptstyle 2};
	(6,7.7)*{\scriptstyle 3};
	(13,7.7)*{\scriptstyle 4};
	(20,7.7)*{\scriptstyle 5};
	(27,7.7)*{\scriptstyle 6};
	(34,7.7)*{\scriptstyle 7};
	(-8,5)*{\scriptstyle\bullet}; 
	(-1,5)*{\scriptstyle\bullet}; 
	(6,5)*{\scriptstyle\bullet} **\dir{};
	(13,5)*{\scriptstyle\bullet} ; **\dir{};
	(20,5)*{\scriptstyle\bullet} **\dir{};
	(27,5)*{\scriptstyle\bullet}; 
	(34,5)*{\scriptstyle\bullet}; 
	(-8,5)*{}="1";
	(-1,5)*{}="2";
	(6,5)*{}="3";
	(13,5)*{}="4";
	(20,5)*{}="5";
	(27,5)*{}="6";
	(34,5)*{}="7";
	(6,2)*{}="3.5";
	"3"; "3.5" **\dir{-};
		(13,2)*{}="4.5";
		"4"; "4.5" **\dir{-};
	(34,2)*{}="7.5";
	"7"; "7.5" **\dir{-};
	"1"; "2" **\crv{(-5.5,2) & (-4,2)};
		"5"; "6" **\crv{(22,1.5) & (25,1.5)};
	(16.5,10)*{}="4.5";
	(16.5,0)*{}="5.5";
	"4.5"; "5.5" **\dir{.};
	\endxy 
~~~~~~~~~~
	\xy
	(-8,7.7)*{\scriptstyle 1};
	(-1,7.7)*{\scriptstyle 2};
	(6,7.7)*{\scriptstyle 3};
	(13,7.7)*{\scriptstyle 4};
	(20,7.7)*{\scriptstyle 5};
	(27,7.7)*{\scriptstyle 6};
	(34,7.7)*{\scriptstyle 7};
	(-8,5)*{\scriptstyle\bullet}; 
	(-1,5)*{\scriptstyle\bullet}; 
	(6,5)*{\scriptstyle\bullet} **\dir{};
	(13,5)*{\scriptstyle\bullet} ; **\dir{};
	(20,5)*{\scriptstyle\bullet} **\dir{};
	(27,5)*{\scriptstyle\bullet}; 
	(34,5)*{\scriptstyle\bullet}; 
	(-8,5)*{}="1";
	(-1,5)*{}="2";
	(6,5)*{}="3";
	(13,5)*{}="4";
	(20,5)*{}="5";
	(27,5)*{}="6";
	(34,5)*{}="7";
	(6,2)*{}="3.5";
	"3"; "3.5" **\dir{-};
		(20,2)*{}="5.5";
		"5"; "5.5" **\dir{-};
	(13,2)*{}="4.5";
	"4"; "4.5" **\dir{-};
	(34,2)*{}="7.5";
	"1"; "2" **\crv{(-5.5,2) & (-4,2)};
	"6"; "7" **\crv{(29,2) & (32,2)};
	(16.5,10)*{}="4.5";
	(16.5,0)*{}="5.5";
	"4.5"; "5.5" **\dir{.};
	\endxy 
	\]
	\[
	\xy
	(-8,7.7)*{\scriptstyle 1};
	(-1,7.7)*{\scriptstyle 2};
	(6,7.7)*{\scriptstyle 3};
	(13,7.7)*{\scriptstyle 4};
	(20,7.7)*{\scriptstyle 5};
	(27,7.7)*{\scriptstyle 6};
	(34,7.7)*{\scriptstyle 7};
	(-8,5)*{\scriptstyle\bullet}; 
	(-1,5)*{\scriptstyle\bullet}; 
	(6,5)*{\scriptstyle\bullet} **\dir{};
	(13,5)*{\scriptstyle\bullet} ; **\dir{};
	(20,5)*{\scriptstyle\bullet} **\dir{};
	(27,5)*{\scriptstyle\bullet}; 
	(34,5)*{\scriptstyle\bullet}; 
	(-8,5)*{}="1";
	(-1,5)*{}="2";
	(6,5)*{}="3";
	(13,5)*{}="4";
	(20,5)*{}="5";
	(27,5)*{}="6";
	(34,5)*{}="7";
	(-8,2)*{}="1.5";
	(6,2)*{}="3.5";
	"1"; "1.5" **\dir{-};
	(13,2)*{}="4.5";
		"4"; "4.5" **\dir{-};
		(20,2)*{}="5.5";
	(27,2)*{}="6.5";
	(34,2)*{}="7.5";
	(34,2)*{}="7.5";
	"7"; "7.5" **\dir{-};
	"5"; "6" **\crv{(22,1.5) & (25,1.5)};
	"2"; "3" **\crv{(1.5,1.5) & (3.5,1.5)};
	(16.5,10)*{}="4.5";
	(16.5,0)*{}="5.5";
	"4.5"; "5.5" **\dir{.};
	\endxy 
~~~~~~~~~~
	\xy
	(-8,7.7)*{\scriptstyle 1};
	(-1,7.7)*{\scriptstyle 2};
	(6,7.7)*{\scriptstyle 3};
	(13,7.7)*{\scriptstyle 4};
	(20,7.7)*{\scriptstyle 5};
	(27,7.7)*{\scriptstyle 6};
	(34,7.7)*{\scriptstyle 7};
	(-8,5)*{\scriptstyle\bullet}; 
	(-1,5)*{\scriptstyle\bullet}; 
	(6,5)*{\scriptstyle\bullet} **\dir{};
	(13,5)*{\scriptstyle\bullet} ; **\dir{};
	(20,5)*{\scriptstyle\bullet} **\dir{};
	(27,5)*{\scriptstyle\bullet}; 
	(34,5)*{\scriptstyle\bullet}; 
	(-8,5)*{}="1";
	(-1,5)*{}="2";
	(6,5)*{}="3";
	(13,5)*{}="4";
	(20,5)*{}="5";
	(27,5)*{}="6";
	(34,5)*{}="7";
	(-8,2)*{}="1.5";
	(6,2)*{}="3.5";
	"1"; "1.5" **\dir{-};
	(13,2)*{}="4.5";
	"4"; "4.5" **\dir{-};
	(20,2)*{}="5.5";
	(27,2)*{}="6.5";
	(34,2)*{}="7.5";
	(34,2)*{}="7.5";
	"2"; "3" **\crv{(1.5,1.5) & (3.5,1.5)};
	(20,2)*{}="5.5";
	"5"; "5.5" **\dir{-};
	(13,2)*{}="4.5";
	"4"; "4.5" **\dir{-};
	"6"; "7" **\crv{(29,2) & (32,2)};
	(16.5,10)*{}="4.5";
	(16.5,0)*{}="5.5";
	"4.5"; "5.5" **\dir{.};
	\endxy 
	\]
	\[
\xy
(-8,7.7)*{\scriptstyle 1};
(-1,7.7)*{\scriptstyle 2};
(6,7.7)*{\scriptstyle 3};
(13,7.7)*{\scriptstyle 4};
(20,7.7)*{\scriptstyle 5};
(27,7.7)*{\scriptstyle 6};
(34,7.7)*{\scriptstyle 7};
(-8,5)*{\scriptstyle\bullet}; 
(-1,5)*{\scriptstyle\bullet}; 
(6,5)*{\scriptstyle\bullet} **\dir{};
(13,5)*{\scriptstyle\bullet} ; **\dir{};
(20,5)*{\scriptstyle\bullet} **\dir{};
(27,5)*{\scriptstyle\bullet}; 
(34,5)*{\scriptstyle\bullet}; 
(-8,5)*{}="1";
(-1,5)*{}="2";
(6,5)*{}="3";
(13,5)*{}="4";
(20,5)*{}="5";
(27,5)*{}="6";
(34,5)*{}="7";
(-8,2)*{}="1.5";
(6,2)*{}="3.5";
"1"; "1.5" **\dir{-};
	(-1,2)*{}="2.5";
	"2"; "2.5" **\dir{-};
(13,2)*{}="4.5";
(20,2)*{}="5.5";
(27,2)*{}="6.5";
(34,2)*{}="7.5";
(34,2)*{}="7.5";
"7"; "7.5" **\dir{-};
"3"; "4" **\crv{(7.5,2) & (11,2)};
"5"; "6" **\crv{(22,1.5) & (25,1.5)};
(16.5,10)*{}="4.5";
(16.5,0)*{}="5.5";
"4.5"; "5.5" **\dir{.};
\endxy 
~~~~~~~~~~
\xy
(-8,7.7)*{\scriptstyle 1};
(-1,7.7)*{\scriptstyle 2};
(6,7.7)*{\scriptstyle 3};
(13,7.7)*{\scriptstyle 4};
(20,7.7)*{\scriptstyle 5};
(27,7.7)*{\scriptstyle 6};
(34,7.7)*{\scriptstyle 7};
(-8,5)*{\scriptstyle\bullet}; 
(-1,5)*{\scriptstyle\bullet}; 
(6,5)*{\scriptstyle\bullet} **\dir{};
(13,5)*{\scriptstyle\bullet} ; **\dir{};
(20,5)*{\scriptstyle\bullet} **\dir{};
(27,5)*{\scriptstyle\bullet}; 
(34,5)*{\scriptstyle\bullet}; 
(-8,5)*{}="1";
(-1,5)*{}="2";
(6,5)*{}="3";
(13,5)*{}="4";
(20,5)*{}="5";
(27,5)*{}="6";
(34,5)*{}="7";
(-8,2)*{}="1.5";
(6,2)*{}="3.5";
"1"; "1.5" **\dir{-};
	(-1,2)*{}="2.5";
	"2"; "2.5" **\dir{-};
(13,2)*{}="4.5";
(20,2)*{}="5.5";
(27,2)*{}="6.5";
(34,2)*{}="7.5";
(34,2)*{}="7.5";
"3"; "4" **\crv{(7.5,2) & (11,2)};
(20,2)*{}="5.5";
"5"; "5.5" **\dir{-};
(13,2)*{}="4.5";
"6"; "7" **\crv{(29,2) & (32,2)};
(16.5,10)*{}="4.5";
(16.5,0)*{}="5.5";
"4.5"; "5.5" **\dir{.};
\endxy 
\]
	
	\caption*{ $\De_{4|3}(0,1,1)\cong \De_{4}(1)\o \De_{3}(1)$}
\end{figure}

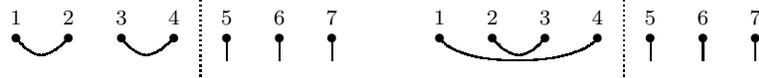
\begin{figure}[h!]
	
	\[
	\xy
	(-8,7.7)*{\scriptstyle 1};
	(-1,7.7)*{\scriptstyle 2};
	(6,7.7)*{\scriptstyle 3};
	(13,7.7)*{\scriptstyle 4};
	(20,7.7)*{\scriptstyle 5};
	(27,7.7)*{\scriptstyle 6};
	(34,7.7)*{\scriptstyle 7};
	(-8,5)*{\scriptstyle\bullet}; 
	(-1,5)*{\scriptstyle\bullet}; 
	(6,5)*{\scriptstyle\bullet} **\dir{};
	(13,5)*{\scriptstyle\bullet} ; **\dir{};
	(20,5)*{\scriptstyle\bullet} **\dir{};
	(27,5)*{\scriptstyle\bullet}; 
	(34,5)*{\scriptstyle\bullet}; 
	(-8,5)*{}="1";
	(-1,5)*{}="2";
	(6,5)*{}="3";
	(13,5)*{}="4";
	(20,5)*{}="5";
	(27,5)*{}="6";
	(34,5)*{}="7";
	(6,2)*{}="3.5";
		(20,2)*{}="5.5";
		"5"; "5.5" **\dir{-};
	(27,2)*{}="6.5";
	"6"; "6.5" **\dir{-};
	(34,2)*{}="7.5";
	"7"; "7.5" **\dir{-};
	"1"; "2" **\crv{(-5.5,2) & (-4,2)};
	"3"; "4" **\crv{(8.5,2) & (10,2)};
	(16.5,10)*{}="4.5";
	(16.5,0)*{}="5.5";
	"4.5"; "5.5" **\dir{.};
	\endxy 
~~~~~~~~~~
	\xy
	(-8,7.7)*{\scriptstyle 1};
	(-1,7.7)*{\scriptstyle 2};
	(6,7.7)*{\scriptstyle 3};
	(13,7.7)*{\scriptstyle 4};
	(20,7.7)*{\scriptstyle 5};
	(27,7.7)*{\scriptstyle 6};
	(34,7.7)*{\scriptstyle 7};
	(-8,5)*{\scriptstyle\bullet}; 
	(-1,5)*{\scriptstyle\bullet}; 
	(6,5)*{\scriptstyle\bullet} **\dir{};
	(13,5)*{\scriptstyle\bullet} ; **\dir{};
	(20,5)*{\scriptstyle\bullet} **\dir{};
	(27,5)*{\scriptstyle\bullet}; 
	(34,5)*{\scriptstyle\bullet}; 
	(-8,5)*{}="1";
	(-1,5)*{}="2";
	(6,5)*{}="3";
	(13,5)*{}="4";
	(20,5)*{}="5";
	(27,5)*{}="6";
	(34,5)*{}="7";
	(-8,2)*{}="1.5";
	(6,2)*{}="3.5";
		(20,2)*{}="5.5";
		"5"; "5.5" **\dir{-};
	(27,2)*{}="6.5";
	"6"; "6.5" **\dir{-};
	(34,2)*{}="7.5";
	"7"; "7.5" **\dir{-};
	"1"; "4" **\crv{(-5.5,1) & (10,1)};
	"2"; "3" **\crv{(1.5,2) & (3.5,2)};
	(16.5,10)*{}="4.5";
	(16.5,0)*{}="5.5";
	"4.5"; "5.5" **\dir{.};
	\endxy 
	\]

	\caption*{ $\De_{4|3}(0,2,0)\cong \De_{4}(2)\o \De_{3}(0)$}
	\caption{ bases of the walled modules}
	\label{fig7}
\end{figure}

\vspace{3cm}

\section*{Acknowledgement}The author acknowledges his supervisor Prof. C.C. Xi.
Also, the author is deeply indebted to Prof. B.M. Deng for introducing the interesting problem with regard to the algebraic structures on Grothendieck groups.

Pei Wang

Basic Courses Department, Beijing Union University, Beijing, 100101, P.R.China

E-mail:wangpei19@mail.bnu.edu.cn


\end{document}